\documentclass[12pt]{amsart}
\usepackage{amssymb,amscd,amsthm,amsmath,color}
\usepackage{fullpage}
\usepackage{epsfig}
\usepackage{graphicx}
\usepackage{xypic}
\usepackage{url}
\usepackage{mathtools}
\usepackage{hyperref}
\usepackage{mathrsfs}
\usepackage{enumerate}

\numberwithin{equation}{section}

\newcommand{\Spec}{\operatorname{Spec}}

\newcommand{\rk}{\operatorname{rk}}

\newcommand{\Eff}{\overline{\operatorname{Eff}}}

\newcommand{\Nef}{\operatorname{Nef}}

\newcommand{\Supp}{\operatorname{Supp}}

\newcommand{\HNP}{\mathrm{HNP}}

\newtheorem{theorem}{Theorem}[section]
\newtheorem{lemma}[theorem]{Lemma}
\newtheorem{proposition}[theorem]{Proposition}
\newtheorem{corollary}[theorem]{Corollary}

\theoremstyle{definition}
\newtheorem{definition}[theorem]{Definition}
\newtheorem{notation}[theorem]{Notation}

\newtheorem{remark}[theorem]{Remark}

\newtheorem{construction}[theorem]{Construction}

\definecolor{ao(english)}{rgb}{0.0, 0.5, 0.0}

\input xy
\xyoption{all}

\makeatletter
\def\dasharrowfill@#1#2#3#4{%
        $\m@th
        \thickmuskip0mu
        \medmuskip\thickmuskip
        \thinmuskip\thickmuskip
        \relax
        #4#1\mkern2mu
        \xleaders\hbox{$#4\mkern2mu#2\mkern2mu$}\hfill
        \mkern2mu
        #3$%
}

\def\dashrightarrowfill@{\dasharrowfill@\relbar\relbar\rightarrow}

\providecommand*\xdashrightarrow[2][]{%
  \ext@arrow 0055{\dashrightarrowfill@}{#1}{#2}}
\makeatother

\begin{document}

\title{Slope stability for normal varieties}

\author{Eric Jovinelly}
\address{Department of Mathematics \\
Brown University \\
Box 1917 \\
151 Thayer Street \\
Providence, RI, 02912}
\email{eric\_jovinelly@brown.edu}

\author{Brian Lehmann}
\address{Department of Mathematics \\
Boston College  \\
Chestnut Hill, MA \, \, 02467}
\email{lehmannb@bc.edu}

\author{Eric Riedl}
\address{Department of Mathematics \\
University of Notre Dame  \\
255 Hurley Hall \\
Notre Dame, IN 46556}
\email{eriedl@nd.edu}

\begin{abstract}
In this expository paper, we review the theory of slope stability of torsion-free sheaves with respect to nef curve classes and extend it to normal varieties.
\end{abstract}

\maketitle

\section{Introduction}

Suppose that $X$ is a smooth projective variety.  In \cite{CP11}, Campana and Peternell introduced the notion of slope stability with respect to a nef curve class on $X$.  Their definition extends the classical notion of Mumford-Takemoto slope stability and is well-suited for applications to birational geometry.  The properties of this construction were studied further in the series of papers \cite{GKP14}, \cite{GKP16} which extend this notion to $\mathbb{Q}$-factorial normal projective varieties.

We generalize this construction one more step: we define slope stability for arbitrary normal projective varieties over fields.  Our motivation is two-fold:
\begin{enumerate}
    \item In characteristic $p$, we currently do not have access to resolution of singularities or $\mathbb{Q}$-factorializations.  Thus it is essential to be able to define slope stability directly for a normal projective variety.
    \item Suppose that $X$ is a smooth projective variety over $\mathbb{C}$.  When the Minimal Model Program terminates with a Fano fibration, there is no guarantee that the fibers of this fibration are $\mathbb{Q}$-factorial. To apply induction on dimension, one would like to work directly with the fibers with no $\mathbb{Q}$-factoriality assumption.
\end{enumerate}
The goal of this paper is to provide a reference for the basic properties of slope stability for normal projective varieties.  While the setting is new, the arguments are not: we emphasize that all the proofs can be found elsewhere in the literature (in slightly different settings).  We have elected to reproduce most of these proofs, giving the paper an expository character.

\subsection{Main idea}
Suppose that $X$ is a normal projective variety of dimension $n$.  When $X$ is not $\mathbb{Q}$-factorial, we must be careful to distinguish between Cartier and Weil divisors.  For a torsion-free coherent sheaf $\mathcal{E}$ on $X$, the first Chern class $c_{1}(\mathcal{E})$ can naturally be constructed as a linear equivalence class of Weil divisors.  Thus to define slope stability, we must use classes that are ``dual'' to Weil divisors (instead of Cartier divisors).

In this paper we will use the space $N^{n-1}(X)_{\mathbb{R}}$ of numerical equivalence classes of weight $(n-1)$ polynomials in Chern classes of vector bundles.  This space naturally carries a nef cone $\Nef^{n-1}(X)$ dual to the pseudo-effective cone of Weil divisors.
For $\alpha \in \Nef^{n-1}(X)$, we define
\begin{equation*}
    \mu_{\alpha}(\mathcal{E}) := \frac{c_{1}(\mathcal{E}) \cdot \alpha}{\rk(\mathcal{E})}.
\end{equation*}
Starting from this definition, the theory develops exactly analogously to the well-known theory of Mumford-Takemoto stability for classes $\alpha = H^{n-1}$ where $H$ is an ample Cartier divisor.  (We emphasize that this setting is strictly more general than working with products of $(n-1)$ nef divisors, even when $X$ is smooth; see Remark \ref{rema:productofnefs}.)

In Section \ref{sect:geomop}, we study the behavior of semistability under certain geometric operations (such as tensor products and finite pullbacks).  The main focus of this section is characteristic $p$, where the theory is more subtle.

\subsection{History}

We give a quick and focused overview of the history of the topics addressed in this paper: the algebraic theory of slope stability for sheaves and for principal bundles.

\cite{Mumford63} introduced the theory of stability for vector bundles on curves in connection with moduli spaces of sheaves.  This theory was further developed and extended to higher dimensions in  \cite{Seshadri67,Takemoto72,Maruyama77,Maruyama78}.  The structural theory of Harder-Narasimhan filtrations was developed in \cite{HN74} and further expanded in \cite{Shatz77}.  The connection between stability and restrictions to complete intersection curves was explored by \cite{MR81b, Flenner84}.  A modern account of stability is presented in \cite{HL97}.

Gradually the theory of slope stability was expanded to new situations.  \cite{Miyaoka87b,Langer04} work with curve classes obtained by intersecting nef divisors and \cite{Langer24,Langer25} extend this approach to normal varieties.  \cite{CP11} works with movable curve classes on smooth projective varieties; the foundations were further developed in \cite{GKP14} and \cite{GKP16}.

In parallel, the theory of slope stability for principal $G$-bundles was first developed in \cite{Ramanathan75} and used to construct moduli spaces in \cite{Ramanathan96a, Ramanthan96b}.  Over time, the theory was extended to higher dimensions and normal varieties in the papers \cite{Langer04b, Langer05, GS05, CH06, GLSS08}.  A Harder-Narasimhan structure theory for principal $G$-bundles was developed in \cite{Behrend95, MS02, BH04, Langer05}.

\bigskip

The properties of semistability in characteristic $p$ are more subtle than in characteristic $0$.  One input into this theory comes from representation theory.   A key tool is GIT, and in particular, the theory of instability flags developed in \cite{Kempf78, Rousseau77}.  The relevance to stability was first described in \cite{RR84}.  Some additional results in this direction are \cite{Langer05,CH06}.  Building on work of \cite{Serre94}, \cite{IMP03} studies semistability of low height representations.  This direction is pursued further in \cite{BP03, BH04}.

A second input into characteristic $p$ semistability is Frobenius techniques.  The important role of the Frobenius map was recognized as early as \cite{Gieseker73}.  Some key ingredients in the modern theory are the Cartier correspondence (as formulated in \cite{Katz70}) and foliation theory as developed in \cite{Miyaoka87, Ekedahl87}.  \cite{Langer04} leveraged Frobenius techniques to study semistable bundles in characteristic $p$; in particular, Langer related the failure of strong semistability to slopes of the cotangent bundle. Some other relevant papers in this direction are \cite{MR83, SB98, JRXY06, Sun10,Trivedi12,Langer15}. For principal bundles, strong semistability was studied in \cite{Sun99,BP04,Langer05}.

\subsection{Overview of paper}

After some preliminary discussion, in Section \ref{sect:slopestability} we prove the existence of Harder-Narasimhan filtrations in our general setting.  The results and arguments are essentially the same as those established in \cite{HL97}, \cite{CP11}, \cite[Appendix A]{GKP14}, and \cite{GKP16} which we follow closely.  In Section \ref{sect:geomop} we turn to the study of Frobenius pullbacks and tensor products in characteristic $p$.  Again, the results and arguments are the same as those established in the papers \cite{RR84}, \cite{IMP03}, \cite{Langer04}, \cite{Langer04b}, \cite{Langer15} and related works.

We work over an arbitrary field in Section \ref{sect:slopestability} and over an algebraically closed field of characteristic $p$ in Section \ref{sect:geomop}.

\subsection{Acknowledgements}

Anthropic's Claude LLM was used to help understand and check the material in Section \ref{sect:tensorprods} and to find references.  Both Claude and ChatGPT were used to proofread the final draft.  All the writing (and any mistakes) are due to the authors.

\section{Background}

Throughout we work over a field $\mathbb{K}$.  A variety will mean an irreducible reduced separated scheme of finite type over $\mathbb{K}$.  A big open subset of a variety is any non-empty open subset whose complement has codimension $\geq 2$.

Suppose that $\mathcal{E}$ is a coherent sheaf on a normal variety $X$.  We denote by $(\mathcal{E})_{tf}$ the torsion-free sheaf obtained by quotienting $\mathcal{E}$ by its torsion subsheaf.  We denote by $\widehat{\mathcal{E}}$ the reflexive hull of $\mathcal{E}$, i.e.~$\widehat{\mathcal{E}} = \mathcal{E}^{\vee \vee}$.  We use similar notation to denote the reflexive hull of the result of a sheaf operation, e.g.~$\mathcal{F} \widehat{\otimes} \mathcal{E} := \widehat{\mathcal{F} \otimes \mathcal{E}}$, $f^{\widehat{*}}\mathcal{E} := \widehat{f^{*}\mathcal{E}}$, and so on.

Suppose $X$ is a variety over an algebraically closed field of characteristic $p$.  We denote by $F_{abs}: X \to X$ the absolute Frobenius and by $F^{e}_{abs}$ its $e$th iterate.  We denote by $F_{rel}: X \to X^{(1)}$ the relative Frobenius over the ground field, i.e.~the morphism induced by the base-change diagram
\begin{equation*}
\xymatrix{
X \ar@/^2pc/[rr]^{F_{abs}} \ar[r]_{F_{rel}} \ar[dr] & X^{(1)}  \ar[r]^{g}  \ar[d] & X \ar[d] \\
& \Spec(\mathbb{K}) \ar[r]_{F_{abs}} & \Spec(\mathbb{K})
}
\end{equation*}
Note that $F_{rel}$ is defined over the ground field.  Since the map $g$ is an isomorphism, for every sheaf $\mathcal{E}$ on $X$ there is a sheaf $\mathcal{E}^{(1)}$ on $X^{(1)}$ such that $F_{abs}^{*}\mathcal{E} \cong F_{rel}^{*}\mathcal{E}^{(1)}$.  We will use this fact implicitly several times.

Since the ground field is perfect, we can iterate this construction to get a sequence of maps
\begin{equation*}
\ldots X^{(-2)} \xrightarrow{F_{rel}}  X^{(-1)} \xrightarrow{F_{rel}} X \xrightarrow{F_{rel}} X^{(1)} \xrightarrow{F_{rel}} X^{(2)} \xrightarrow{F_{rel}} \ldots
\end{equation*}

\section{Numerical spaces and cones}

Let $X$ be a projective variety of dimension $n$ over the field $\mathbb{K}$.  We let $A_{k}(X)$ denote the Chow group of $k$-cycles on $X$ up to rational equivalence.

We let $N^{1}(X)_{\mathbb{Z}}$ denote the additive group of Cartier divisors up to numerical equivalence.  We let $N_{1}(X)_{\mathbb{Z}}$ denote the additive group of $1$-cycles up to numerical equivalence.  By tensoring these groups with $\mathbb{R}$ we get finite dimensional vector spaces $N^{1}(X)_{\mathbb{R}}$ and $N_{1}(X)_{\mathbb{R}}$ which are dual under the intersection pairing.  We denote by $N^{1}(X)_{\mathbb{Q}}, N_{1}(X)_{\mathbb{Q}}$ the $\mathbb{Q}$-analogues.

We let $N_{n-1}(X)_{\mathbb{Z}}$ denote the additive group of Weil divisors up to numerical equivalence (that is, up to intersecting against homogeneous polynomials of weighted degree $(n-1)$ in Chern classes of vector bundles on $X$; see \cite{FL17}).  We let $N^{n-1}(X)_{\mathbb{Z}}$ denote the additive group of homogeneous polynomials in Chern classes of vector bundles of weighted degree $(n-1)$ up to numerical equivalence.  By tensoring these groups with $\mathbb{R}$ we get finite-dimensional vector spaces $N_{n-1}(X)_{\mathbb{R}}$ and $N^{n-1}(X)_{\mathbb{R}}$ which are dual under the intersection pairing.  We denote by $N_{n-1}(X)_{\mathbb{Q}}, N^{n-1}(X)_{\mathbb{Q}}$ the $\mathbb{Q}$-analogues.

\begin{remark} \label{rema:insmoothlocus}
Suppose that $C$ is a regular projective curve and $s: C \to X$ is a map whose image lies in the regular locus of $X$.  Note that $s$ is an lci morphism of codimension $n-1$.  Then \cite[Example 19.2.3]{Fulton84} shows that the Gysin map $s^{*}: A_{n-1}(X) \to A_{0}(C)$ descends to numerical groups.  In this way, any curve $C$ in the regular locus of $X$ naturally defines a class $s_{*}C \in N^{n-1}(X)_{\mathbb{Q}}$.
\end{remark}

As explained in \cite[Example 2.7]{FL17} the Cartier-to-Weil map obtained by capping against $[X]$ induces an injective linear map $\cap [X]: N^{1}(X)_{\mathbb{R}} \to N_{n-1}(X)_{\mathbb{R}}$, or dually, a surjective linear map $\cap [X]: N^{n-1}(X)_{\mathbb{R}} \to N_{1}(X)_{\mathbb{R}}$.  When $X$ is a $\mathbb{Q}$-factorial projective variety both these maps are isomorphisms.

Given a morphism $f: X \to Y$ of projective varieties, we get maps $f_{*}: N_{k}(X)_{\mathbb{Z}} \to N_{k}(Y)_{\mathbb{Z}}$ induced by pushforward of cycles and maps $f^{*}: N^{k}(Y)_{\mathbb{Z}} \to N^{k}(X)_{\mathbb{Z}}$ induced by pullback of vector bundles.  We also denote by $f_{*}$ and $f^{*}$ the induced maps after tensoring by $\mathbb{R}$; then $f_{*}$ and $f^{*}$ are dual maps of vector spaces.

\subsection{Cones of divisors and curves}
We next define several ``positive'' cones of divisors and curves.  We will define these cones in $N_{k}(X)_{\mathbb{R}}$ and in $N^{k}(X)_{\mathbb{R}}$ for $k = 1, n-1$.

\begin{definition}
Let $X$ be a projective variety.  For $k=1,n-1$, the pseudo-effective cone of $k$-cycles $\Eff_{k}(X) \subset N_{k}(X)_{\mathbb{R}}$ is the closure of the cone generated by classes of effective $k$-cycles.

These cones are full-dimensional, closed, convex, and pointed in $N_{k}(X)_{\mathbb{R}}$ by \cite[Theorem 1.4.(1)]{FL17}. (While the cited reference assumes the field is algebraically closed, the argument works equally well over all fields.)
\end{definition}

We denote by $\Eff_{k}(X)_{\mathbb{Z}}$ or $\Eff_{k}(X)_{\mathbb{Q}}$ the intersection of $\Eff_{k}(X)$ with $N_{k}(X)_{\mathbb{Z}}$ or $N_{k}(X)_{\mathbb{Q}}$.  Given a morphism $f: X \to Y$ of projective varieties, the pushforward map $f_{*}$ satisfies $f_{*}\Eff_{k}(X) \subset  \Eff_{k}(Y)$.  If $f$ is surjective, then these two sets are equal by \cite[Theorem 1.4.(4)]{FL17}.

\begin{definition}
Let $X$ be a projective variety.  For $k=1,n-1$ the nef cone $\Nef^{k}(X) \subset N^{k}(X)_{\mathbb{R}}$ is the dual of the pseudo-effective cone $\Eff_{k}(X)$. 

These cones are full-dimensional, closed, convex, and pointed in $N^{k}(X)_{\mathbb{R}}$.
\end{definition}

We denote by $\Nef^{k}(X)_{\mathbb{Z}}$ or $\Nef^{k}(X)_{\mathbb{Q}}$ the intersection of $\Nef^{k}(X)$ with $N^{k}(X)_{\mathbb{Z}}$ or $N^{k}(X)_{\mathbb{Q}}$.  Given a morphism $f: X \to Y$ of projective varieties, the pullback map $f^{*}$ satisfies $f^{*}\Nef^{k}(Y) \subset \Nef^{k}(X)$.  If $f$ is surjective, then $f^{*}\Nef^{k}(Y) = \Nef^{k}(X) \cap f^{*}N^{k}(Y)_{\mathbb{R}}$.

\begin{remark}
    As in Remark \ref{rema:insmoothlocus}, a curve $C$ in the regular locus of $X$ that deforms to dominate $X$ naturally defines an element of $\Nef^{n-1}(X)_{\mathbb{Q}}$.
\end{remark}

\begin{remark} \label{rema:productofnefs}
    Given any $(n-1)$ nef $\mathbb{R}$-Cartier divisors $A_{1},\ldots,A_{n-1}$, the product $A_{1} \cdot \ldots \cdot A_{n-1}$ is an element of $\Nef^{n-1}(X)$.  However, the cone $\Nef^{n-1}(X)$ is often larger than the cone of complete intersections of nef classes (which need not even be convex).  This is true even when $X$ is smooth; see \cite[Example 9.1]{LX16} for an explicit toric example.
\end{remark}

\section{Slope stability} \label{sect:slopestability}

In this section we work over an arbitrary field $\mathbb{K}$.

\begin{construction}[First Chern class]
Let $X$ be a normal projective variety of dimension $n$.  For a torsion-free coherent sheaf $\mathcal{F}$ on $X$ we define its first Chern class $c_{1}(\mathcal{F}) \in A_{n-1}(X)$ in the following way.  There is a big open subset $U \subset X$ such that $\mathcal{F}|_{U}$ is locally free.  Then $\bigwedge^{\rk \mathcal{F}} \mathcal{F}|_{U}$ is a line bundle $\mathcal{L}$ on $U$.  Since $A_{n-1}(U) \cong A_{n-1}(X)$, the line bundle $\mathcal{L}$ uniquely defines a rational equivalence class of Weil divisors on $X$ (obtained by taking closures).  It is clear that the resulting class is independent of the choice of big open $U$; we call this class $c_{1}(\mathcal{F}) \in A_{n-1}(X)$.  We define $\det(\mathcal{F})$ to be the rank $1$ reflexive sheaf corresponding to this Weil divisor class.
\end{construction}

\begin{remark} \label{rema:bigopen}
    Suppose that $\mathcal{E}$ is a torsion-free sheaf.  Then $c_{1}(\mathcal{E})$ is determined by the restriction of $\mathcal{E}$ to any big open subset $U \subset X$.  In particular, $c_{1}(\mathcal{E}) = c_{1}(\widehat{\mathcal{E}})$.  We will frequently use this observation without mention in the future.
\end{remark}

The following proposition verifies a few standard properties of the first Chern class in this general setting.

\begin{proposition} \label{prop:c1properties}
Let $X$ be a normal projective variety of dimension $n$.  Then:
\begin{enumerate}
\item $c_{1}$ is additive in exact sequences of torsion-free sheaves.
\item $c_{1}$ of a globally generated torsion-free sheaf is the class of an effective Weil divisor.
\item If $\mathcal{F} \subset \mathcal{F}'$ is an inclusion of torsion-free sheaves of equal rank, then $c_{1}(\mathcal{F}') - c_{1}(\mathcal{F})$ is the class of an effective Weil divisor whose support is the codimension $1$ part of $\Supp(\mathcal{F}'/\mathcal{F})$.
\item If $\mathcal{G}_{1},\mathcal{G}_{2}$ are torsion-free subsheaves of a torsion-free sheaf $\mathcal{F}$, then
\begin{equation*}
c_{1}(\mathcal{G}_{1} + \mathcal{G}_{2}) = c_{1}(\mathcal{G}_{1}) + c_{1}(\mathcal{G}_{2}) - c_{1}(\mathcal{G}_{1} \cap \mathcal{G}_{2}).
\end{equation*}
\end{enumerate}
\end{proposition}

\begin{proof}
Let $U \subset X$ denote the big open subset obtained by intersecting the regular locus with the locus where all sheaves of interest are locally free.  Then there is an isomorphism $A_{n-1}(U) \cong A_{n-1}(X)$ and the effectiveness of a rational equivalence class is preserved by this identification.  All the statements can be verified in the usual way on the regular locus $U$ and then extended to $X$.  
\end{proof}

Let $X$ be a normal projective variety and let $\alpha \in \Nef^{n-1}(X)$.  For any torsion-free coherent sheaf $\mathcal{E}$ we define the slope of $\mathcal{E}$ with respect to $\alpha$ to be
\begin{equation*}
\mu_{\alpha}(\mathcal{E}) = \frac{c_{1}(\mathcal{E}) \cdot \alpha}{\rk(\mathcal{E})}
\end{equation*}

\begin{lemma} \label{lemm:exactsequence}
    Let $X$ be a normal projective variety and let $\alpha \in \Nef^{n-1}(X)$.  Suppose we have an exact sequence of non-zero torsion-free sheaves
    \begin{equation*}
        0 \to \mathcal{F} \to \mathcal{E} \to \mathcal{Q} \to 0.
    \end{equation*}
    Then the slope of $\mathcal{E}$ lies between the slopes of $\mathcal{F}$ and $\mathcal{Q}$.
\end{lemma}

\begin{proof}
    The slope of $\mathcal{E}$ is the rank-weighted average of the slopes of $\mathcal{F},\mathcal{Q}$ and the claim follows.
\end{proof}

\begin{definition}
Let $X$ be a normal projective variety, $\mathcal{E}$ a non-zero torsion-free coherent sheaf on $X$, and $\alpha \in \Nef^{n-1}(X)$.  We say that $\mathcal{E}$ is $\alpha$-stable if
\begin{equation*}
\mu_{\alpha}(\mathcal{F}) < \mu_{\alpha}(\mathcal{E})
\end{equation*}
for every non-zero torsion-free coherent subsheaf $\mathcal{F} \subset \mathcal{E}$ with $\rk(\mathcal{F}) < \rk(\mathcal{E})$.  
We say that $\mathcal{E}$ is $\alpha$-semistable if
\begin{equation*}
\mu_{\alpha}(\mathcal{F}) \leq \mu_{\alpha}(\mathcal{E})
\end{equation*}
for every non-zero torsion-free coherent subsheaf $\mathcal{F} \subset \mathcal{E}$.  (By Proposition \ref{prop:c1properties}.(3) it suffices to check the subsheaves with $\rk(\mathcal{F}) < \rk(\mathcal{E})$.)  We will drop the subscript and prefix ``$\alpha$'' when it is clear from the context.
\end{definition}

\begin{lemma} \label{lemm:semistabpreserved}
    Let $X$ be a normal projective variety and let $\alpha \in \Nef^{n-1}(X)$.
    \begin{enumerate}
        \item If $\mathcal{E}_{1},\mathcal{E}_{2}$ are semistable sheaves of the same slope $\mu$, then $\mathcal{E}_{1} \oplus \mathcal{E}_{2}$ is again semistable of slope $\mu$.
        \item Every rank $1$ torsion-free sheaf on $X$ is semistable.
    \end{enumerate}
\end{lemma}

\begin{proof}
    (1) Lemma \ref{lemm:exactsequence} shows that $\mu_{\alpha}(\mathcal{E}_{1} \oplus \mathcal{E}_{2}) = \mu$.  Suppose $\mathcal{F} \subset \mathcal{E}_{1} \oplus \mathcal{E}_{2}$ is a non-zero subsheaf.  Letting $p_{2}$ denote the second projection map, we have an exact sequence
    \begin{equation*}
        0 \to \mathcal{F} \cap \mathcal{E}_{1} \to \mathcal{F} \to p_{2}(\mathcal{F}) \to 0.
    \end{equation*}
    If one of the end terms is zero, then $\mathcal{F}$ is a subsheaf of $\mathcal{E}_{1}$ or $\mathcal{E}_{2}$ and we conclude that $\mu_{\alpha}(\mathcal{F}) \leq \mu$ by semistability of $\mathcal{E}_{1},\mathcal{E}_{2}$.  Otherwise, the semistability of $\mathcal{E}_{1},\mathcal{E}_{2}$ implies that both of the end terms have slope at most $\mu$.  We conclude that $\mathcal{F}$ also has slope at most $\mu$ by Lemma \ref{lemm:exactsequence}.
    
    (2) Immediate from the definition.
\end{proof}

\subsection{Maximal destabilizing subsheaves}

\begin{proposition}[{\cite[Proposition A.2]{GKP14}}] \label{prop:maxslopeexists}
Let $X$ be a normal projective variety and let $\alpha \in \Nef^{n-1}(X)$.  For any non-zero torsion-free sheaf $\mathcal{E}$, the quantity
\begin{equation*}
\mu^{max}_{\alpha}(\mathcal{E}) := \sup\{ \mu_{\alpha}(\mathcal{G}) \, | \, \mathcal{G} \subset \mathcal{E}, \, \mathcal{G} \textrm{ torsion-free and non-zero} \}
\end{equation*}
is finite and is achieved by a subsheaf $\mathcal{F}$.
\end{proposition}

The following proof is essentially the same as the one given in \cite[Appendix A]{GKP14} or in \cite[Proposition 2.23]{GKP16}.

\begin{proof}
We first prove that the slopes of torsion-free subsheaves of $\mathcal{E}$ are bounded above.  Fix an ample Cartier divisor $H$ on $X$.  For some choice of positive integers $m,k$ there is an inclusion $\mathcal{E} \subset \mathcal{O}_{X}(mH)^{\oplus k}$.  By Lemma \ref{lemm:semistabpreserved} we see that for every non-zero torsion-free subsheaf $\mathcal{G} \subset \mathcal{E}$ we have $c_{1}(\mathcal{G}) \cdot \alpha \leq kmH \cdot \alpha$, proving the desired upper bound on the slopes of subsheaves of $\mathcal{E}$.

Assume for a contradiction that there is no torsion-free subsheaf achieving the maximal slope.  Choose an infinite sequence of non-zero torsion-free subsheaves $\{ \mathcal{G}_{i} \}$ whose slopes converge to $\mu^{max}_{\alpha}(\mathcal{E})$.  By passing to a subsequence, we may suppose that all the $\mathcal{G}_{i}$ have the same rank $r$.  We may also ensure that $r$ is the maximal possible rank of any constant rank sequence whose slope converges to $\mu^{max}_{\alpha}(\mathcal{E})$.  We may also suppose that each $\mathcal{G}_{i}$ is saturated (since by Proposition \ref{prop:c1properties}.(3) taking saturations can only increase the slope).  Again passing to a subsequence, we may suppose that the slopes of the $\mathcal{G}_{i}$ are strictly increasing.  Finally, we may suppose that the slope of $\mathcal{G}_{i}$ satisfies $\mu(\mathcal{G}_{i}) > \mu^{max}_{\alpha}(\mathcal{E}) - \frac{1}{i}$.

We claim that $\mathcal{G}_{i} \not \subset \mathcal{G}_{j}$ for $i<j$.  Indeed, if we had a containment, then we would also have $\mathcal{G}_{j}/\mathcal{G}_{i} \subset \mathcal{E}/\mathcal{G}_{i}$.  However, the first sheaf is torsion (since it is a quotient of two sheaves of the same rank) and the second sheaf is torsion-free by construction.  Thus we must have $\mathcal{G}_{i} = \mathcal{G}_{j}$, contradicting the fact that the slopes are strictly increasing.  This proves the claim.

Now define $\mathcal{H}_{i,j} := \mathcal{G}_{i} + \mathcal{G}_{j} \subset \mathcal{E}$.  Since $\mathcal{G}_{j}$ is saturated in $\mathcal{E}$, the fact that $\mathcal{G}_{i} \not \subset \mathcal{G}_{j}$ implies that $\mathcal{H}_{i,j}$ has larger rank than both $\mathcal{G}_{i}$ and $\mathcal{G}_{j}$.  Furthermore by Proposition \ref{prop:c1properties}.(4) we have
\begin{align*}
    \rk(\mathcal{H}_{i,j}) \cdot \mu_{\alpha}(\mathcal{H}_{i,j}) & = r(\mu_{\alpha}(\mathcal{G}_{i}) + \mu_{\alpha}(\mathcal{G}_{j})) - \rk(\mathcal{G}_{i} \cap \mathcal{G}_{j}) \mu_{\alpha}(\mathcal{G}_{i} \cap \mathcal{G}_{j}) \\
    & \geq r\left(2\mu_{\alpha}^{max}(\mathcal{E}) - \tfrac{1}{i} - \tfrac{1}{j} \right)  - \rk(\mathcal{G}_{i} \cap \mathcal{G}_{j}) \mu_{\alpha}^{max}(\mathcal{E}) \\
    & = \rk(\mathcal{H}_{i,j})\mu_{\alpha}^{max}(\mathcal{E}) - r \left( \tfrac{1}{i} + \tfrac{1}{j} \right).
\end{align*}
Thus for the sequence of sheaves $\{\mathcal{H}_{j,j+1}\}$ the slopes converge to $\mu^{max}_{\alpha}(\mathcal{E})$ and each term has rank strictly larger than $r$.  By passing to a subsequence, we may get a sequence of constant rank.  This contradicts the choice of $r$ as the maximal rank amongst all such sequences, finishing the proof.
\end{proof}

The following results are immediate.

\begin{corollary} \label{coro:maxdestabss}
    Let $X$ be a normal projective variety and let $\alpha \in \Nef^{n-1}(X)$.  Fix a non-zero torsion-free sheaf $\mathcal{E}$.  Every non-zero torsion-free subsheaf $\mathcal{F} \subset \mathcal{E}$ with $\mu_{\alpha}(\mathcal{F}) = \mu_{\alpha}^{max}(\mathcal{E})$ is semistable.
\end{corollary}

\begin{corollary} \label{coro:mapandslope}
    Let $X$ be a normal projective variety and let $\alpha \in \Nef^{n-1}(X)$.  Suppose that $\mathcal{E}$ is a semistable sheaf.  If there is a non-zero homomorphism $\mathcal{E} \to \mathcal{G}$ to a non-zero torsion-free sheaf $\mathcal{G}$ then $\mu(\mathcal{E}) \leq \mu^{max}(\mathcal{G})$.
\end{corollary}

We next verify that there is a ``largest'' subsheaf $\mathcal{E}_{max}$ achieving the maximum of Proposition \ref{prop:maxslopeexists}.

\begin{proposition}[{\cite[Corollary 2.25]{GKP16}}] \label{prop:maxdestab}
    Let $X$ be a normal projective variety and let $\alpha \in \Nef^{n-1}(X)$.  Fix a non-zero torsion-free sheaf $\mathcal{E}$.  There is a unique saturated subsheaf $\mathcal{E}_{max} \subset \mathcal{E}$ such that
    \begin{itemize}
        \item $\mu(\mathcal{E}_{max}) = \mu^{max}(\mathcal{E})$, and
        \item every non-zero subsheaf $\mathcal{F} \subset \mathcal{E}$ with $\mu(\mathcal{F}) = \mu^{max}(\mathcal{E})$ satisfies $\mathcal{F} \subset \mathcal{E}_{max}$.
    \end{itemize}
\end{proposition}

\begin{proof}
    Set $\mu = \mu^{max}(\mathcal{E})$.  Let $\mathcal{E}_{max} \subset \mathcal{E}$ be the saturation of any non-zero torsion-free subsheaf of maximal rank amongst all subsheaves of slope $\mu$.  We claim that this subsheaf has the desired property.  Indeed, if $\mathcal{F} \subset \mathcal{E}$ also has slope $\mu$ but is not contained in $\mathcal{E}_{max}$, then by Proposition \ref{prop:c1properties}.(4),  $\mathcal{F} + \mathcal{E}_{max}$ would have the same slope as $\mathcal{E}_{max}$ but strictly larger rank, a contradiction.  The uniqueness is clear from the defining property of $\mathcal{E}_{max}$.
\end{proof}

\begin{definition}
    In the setting of Proposition \ref{prop:maxdestab}, the subsheaf $\mathcal{E}_{max}$ is called the maximal destabilizing subsheaf of $\mathcal{E}$ with respect to $\alpha$.  Note that it is semistable by Corollary \ref{coro:maxdestabss}.
\end{definition}

\subsection{Harder-Narasimhan filtrations}
We next construct the Harder-Narasimhan filtration of a torsion-free sheaf $\mathcal{E}$.

\begin{definition}
    Let $\mathcal{E}$ be a torsion-free sheaf and let $\alpha \in \Nef^{n-1}(X)$.  An $\alpha$-Harder-Narasimhan filtration for $\mathcal{E}$ is a filtration
    \begin{equation*}
        0 = \mathcal{F}_{0} \subset \mathcal{F}_{1} \subset \ldots \subset \mathcal{F}_{k} = \mathcal{E}
    \end{equation*}
    satisfying the following properties:
    \begin{enumerate}
        \item Each subsheaf $\mathcal{F}_{i}$ is saturated in $\mathcal{E}$.
        \item Each successive quotient $\mathcal{F}_{i}/\mathcal{F}_{i-1}$ is torsion-free and semistable.
        \item The sequence of numbers $\mu_{\alpha}(\mathcal{F}_{i}/\mathcal{F}_{i-1})$ is strictly decreasing as $i$ increases.
    \end{enumerate}
    We will drop the prefix ``$\alpha$'' when it is clear from the context.
\end{definition}

\begin{theorem}
    Let $X$ be a normal projective variety and let $\alpha \in \Nef^{n-1}(X)$.  Then every non-zero torsion-free sheaf $\mathcal{E}$ admits a unique Harder-Narasimhan filtration.
\end{theorem}

\begin{proof}
To prove the existence, we inductively define $\mathcal{F}_{i} \subset \mathcal{E}$ to be the preimage of the maximal destabilizing subsheaf in $\mathcal{E}/\mathcal{F}_{i-1}$.  The semistability of $\mathcal{F}_{i}/\mathcal{F}_{i-1}$ was proved in Corollary \ref{coro:maxdestabss}.  If $\mu(\mathcal{F}_{i+1}/\mathcal{F}_{i})$ were not strictly smaller than $\mu(\mathcal{F}_{i}/\mathcal{F}_{i-1})$, then by Lemma \ref{lemm:exactsequence} $\mathcal{F}_{i+1}/\mathcal{F}_{i-1}$ would have slope at least as large as $\mathcal{F}_{i}/\mathcal{F}_{i-1}$.  This would contradict the defining feature of maximal destabilizing subsheaves identified in Proposition \ref{prop:maxdestab}.

The uniqueness of the sequence follows from an easy induction using Corollary \ref{coro:mapandslope}.
\end{proof}

\begin{definition}
Let $X$ be a normal projective variety and let $\mathcal{E}$ be a torsion-free sheaf on $X$.  For any class $\alpha \in \Nef^{n-1}(X)$, we define $\mu^{min}_{\alpha}(\mathcal{E})$ to be the minimal slope of any non-zero torsion-free quotient of $\mathcal{E}$ (or equivalently, the slope of the quotient of the last two terms in the Harder-Narasimhan filtration).
\end{definition}

We next verify a couple basic properties of Harder-Narasimhan filtrations.

\begin{lemma} \label{lemm:slopeinequalityzerohom}
    Let $X$ be a normal projective variety and let $\alpha \in \Nef^{n-1}(X)$ be a nef class.  Suppose that $\mathcal{E}_{1},\mathcal{E}_{2}$ are torsion-free sheaves on $X$ satisfying $\mu^{min}(\mathcal{E}_{1}) > \mu^{max}(\mathcal{E}_{2})$.  Then $\mathrm{Hom}(\mathcal{E}_{1},\mathcal{E}_{2}) = 0$.
\end{lemma}

\begin{proof}
    Suppose for a contradiction that there were a non-zero homomorphism $\phi: \mathcal{E}_{1} \to \mathcal{E}_{2}$.  Then the image $\mathcal{G}$ is a non-zero torsion-free sheaf.  Since we have a surjection $\mathcal{E}_{1} \to \mathcal{G}$, we conclude $\mu^{min}(\mathcal{E}_{1}) \leq \mu(\mathcal{G})$.  Since we have an injection $\mathcal{G} \to \mathcal{E}_{2}$, we conclude that $\mu(\mathcal{G}) \leq \mu^{max}(\mathcal{E}_{2})$.  These two inequalities provide the contradiction.
\end{proof}

\begin{lemma} \label{lemm:dual}
    Let $X$ be a normal projective variety and let $\alpha \in \Nef^{n-1}(X)$ be a nef class.  Suppose that $\mathcal{E}$ is a torsion-free sheaf on $X$ with Harder-Narasimhan filtration
    \begin{equation*}
        0 = \mathcal{F}_{0} \subset \mathcal{F}_{1} \subset \ldots \subset \mathcal{F}_{k} = \mathcal{E}.
    \end{equation*}
    Then the $i$th term in the Harder-Narasimhan filtration of $\mathcal{E}^{\vee}$ is $(\mathcal{E}/\mathcal{F}_{k-i})^{\vee}$.
\end{lemma}

\begin{proof}
    Suppose that $\mathcal{G} \subset \mathcal{E}^{\vee}$ is a saturated subsheaf.  Then we have a generically surjective morphism $\phi: \mathcal{E} \to \mathcal{G}^{\vee}$ whose cokernel is supported on a closed set of codimension $\geq 2$. Note that $\mu(\textrm{im} \, \phi) = \mu(\mathcal{G}^{\vee})$.  Conversely, given any torsion-free quotient $\mathcal{E} \to \mathcal{Q}$, we get an injection $\mathcal{Q}^{\vee} \to \mathcal{E}^{\vee}$.  The image is a saturated subsheaf of $\mathcal{E}^{\vee}$.
    
    Recall that for any torsion-free sheaf $\mathcal{G}$ we have $c_{1}(\mathcal{G}) = -c_{1}(\mathcal{G}^{\vee})$; indeed, this computation can be executed on the big open subset where $\mathcal{G}$ is locally free.  Combined with the previous paragraph, we get a precise correspondence between slopes of torsion-free quotients of $\mathcal{E}$ and saturated subsheaves of $\mathcal{E}^{\vee}$, proving the result.
\end{proof}

\subsection{Birational properties}

The following lemma summarizes the birational properties of the slope function.

\begin{lemma}[{\cite[Lemma A.19, Lemma A.20, Lemma A.21]{GKP14}}] \label{lemm:slopebirational}
Let $\phi: X \to Y$ be a birational morphism of normal projective varieties of dimension $n$ and let $\alpha \in \Nef^{n-1}(Y)$.
\begin{enumerate}
\item For a torsion-free coherent sheaf $\mathcal{F}$ on $X$, we have that $\phi_{*}\mathcal{F}$ is torsion-free and $\mu_{\alpha}(\phi_{*}\mathcal{F}) = \mu_{\phi^{*}\alpha}(\mathcal{F})$.
\item For a torsion-free coherent sheaf $\mathcal{G}$ on $Y$, we have that $(\phi^{*}\mathcal{G})_{tf}$ is torsion-free and $\mu_{\alpha}(\mathcal{G}) = \mu_{\phi^{*}\alpha}((\phi^{*}\mathcal{G})_{tf})$.
\item  Let $U \subset X$ be the open set where $\phi$ is an isomorphism.  Suppose that $\mathcal{E}$ and $\mathcal{E}'$ are two torsion-free sheaves on $X$ such that there is an isomorphism $\mathcal{E}|_{U} \cong \mathcal{E}'|_{U}$.  Then $\mu_{\phi^{*}\alpha}(\mathcal{E}) = \mu_{\phi^{*}\alpha}(\mathcal{E}')$.
\end{enumerate}
\end{lemma}

\begin{proof}
(1) The torsion-freeness of $\phi_{*}\mathcal{F}$ is clear.  To prove the second statement, consider the rank $1$ torsion-free sheaves $\phi_{*}(\det \mathcal{F})$ and $\det(\phi_{*}\mathcal{F})$.  These are isomorphic away from the codimension $2$ set consisting of the $\phi$-exceptional centers, and thus define the same class in $A_{n-1}(Y)$.

(2) The torsion-freeness of $(\phi^{*}\mathcal{G})_{tf}$ is clear.  To prove the second statement, note that $\det(\phi_{*}(\phi^{*}\mathcal{G})_{tf})$ is isomorphic to $\det(\mathcal{G})$ away from the codimension $2$ locus consisting of the locus where $\phi^{-1}$ is not an isomorphism.  Thus $\det(\phi_{*}(\phi^{*}\mathcal{G})_{tf})$ and $\det(\mathcal{G})$ define the same class in $A_{n-1}(Y)$.  We then apply (1).

(3) It is clear that $\det(\mathcal{E})|_{U}$ and $\det(\mathcal{E}')|_{U}$ are isomorphic on $U$.  Using the exact sequence
\begin{equation*}
A_{n-1}(X \backslash U) \to A_{n-1}(X) \to A_{n-1}(U) \to 0
\end{equation*}
we see that the first Chern classes of $\mathcal{E}$ and $\mathcal{E}'$ only differ by the class of a $\phi$-exceptional divisor.  Thus the intersections against $\phi^{*}\alpha$ are the same.
\end{proof}

As a consequence, we have:

\begin{proposition}[\cite{GKP14} Corollary A.22] \label{prop:semistablebirational}
Let $\phi: X \to Y$ be a birational morphism of normal projective varieties of dimension $n$.  Fix a class $\alpha \in \Nef^{n-1}(Y)$.  Let $\mathcal{F}$ be a torsion-free sheaf on $Y$ and let $\mathcal{G}$ be a torsion-free sheaf on $X$ such that $(\phi^{*}\mathcal{F})_{tf}$ and $\mathcal{G}$ are isomorphic on the locus where $\phi$ is an isomorphism.

Given numbers $(r,q) \in \mathbb{N} \times \mathbb{R}$, $\mathcal{F}$ contains a torsion-free subsheaf of rank $r$ and $\alpha$-slope $q$ if and only if $\mathcal{G}$ contains a torsion-free subsheaf of rank $r$ and $\phi^{*}\alpha$-slope $q$.
\end{proposition}

In fact, the proof of \cite[Corollary A.22]{GKP14} shows that we can construct the desired subsheaves of $\mathcal{F}$ and $\mathcal{G}$ by restricting to the locus where $\phi$ is an isomorphism and taking saturated extensions to the ambient projective variety.  This implies that the Harder-Narasimhan filtration has a natural compatibility with birational maps.

\begin{corollary} \label{coro:hnfiltrationbirational}
Let $\phi: X \to Y$ be a birational morphism of normal projective varieties of dimension $n$ and let $\alpha \in \Nef^{n-1}(Y)$.  Let $\mathcal{E}$ be a torsion-free coherent sheaf on $Y$ and let
\begin{equation*}
0 = \mathcal{F}_{0} \subset \mathcal{F}_{1} \subset \ldots \subset \mathcal{F}_{s} = \mathcal{E}
\end{equation*}
be the $\alpha$-Harder-Narasimhan filtration for $\mathcal{E}$.  Let $\mathcal{G}$ be any torsion-free sheaf on $X$ such that $(\phi^{*}\mathcal{E})_{tf}$ and $\mathcal{G}$ are isomorphic on the locus where $\phi$ is an isomorphism.  Then the $\phi^{*}\alpha$-Harder-Narasimhan filtration for $\mathcal{G}$ has the form
\begin{equation*}
0 = \mathcal{F}'_{0} \subset \mathcal{F}'_{1} \subset \ldots \subset \mathcal{F}'_{s} = \mathcal{G}
\end{equation*}
where $\mathcal{F}'_{i}$ is isomorphic to $\mathcal{F}_{i}$ on the locus where $\phi$ is an isomorphism and
\begin{equation*}
\mu_{\alpha}(\mathcal{F}_{i}/\mathcal{F}_{i-1}) = \mu_{\phi^{*}\alpha}(\mathcal{F}'_{i}/\mathcal{F}'_{i-1}).
\end{equation*}
\end{corollary}

\subsection{Continuity}
\cite[Lemma 3.3.3]{Neumann10} proves a fundamental continuity property of the Harder-Narasimhan filtration.  We verify that the same statement holds in our setting.  For the reader's convenience we will give the details of the proof, following \cite{Neumann10} closely.

\begin{theorem}[{\cite[Lemma 3.3.3]{Neumann10}}] \label{theo:finmanydestabilizing}
Let $X$ be a normal projective variety and let $\mathcal{E}$ be a torsion-free sheaf on $X$.  Fix a finite set of nef curve classes $\alpha_{1},\ldots,\alpha_{k} \in \Nef^{n-1}(X)_{\mathbb{Q}}$.  Then as we vary $\beta \in \langle \alpha_{1}, \ldots, \alpha_{k} \rangle$ there are only finitely many $\beta$-maximal destabilizing subsheaves of $\mathcal{E}$.
\end{theorem}

\begin{proof}
Without loss of generality we may assume that $\langle \alpha_{1}, \ldots, \alpha_{k} \rangle$ is a full-dimensional simplicial subcone $\mathcal{C}$ of $\Nef^{n-1}(X)$.  Let $\{D_{1},\ldots,D_{k} \} \subset N_{n-1}(X)_{\mathbb{Q}}$ be the basis dual to the $\alpha_{i}$.  Note that there is some fixed integer $m$ such that every class in $N_{n-1}(X)_{\mathbb{Z}}$ can be written as a linear combination of the $D_{i}$ with coefficients in $\frac{1}{m}\mathbb{Z}$.

Since the maximal destabilizing subsheaf is invariant under rescaling of the nef curve class, it suffices to consider $\beta$ in a fixed compact slice of the nef cone.   Suppose for a contradiction that there is an infinite sequence $\{ \beta_{j} \}$ of distinct classes in our fixed compact slice such that there are infinitely many distinct $\beta_{j}$-maximal destabilizing subsheaves of $\mathcal{E}$.  We denote the subsheaf for $\beta_{j}$ by $\mathcal{G}_{j}$.  Passing to a subsequence, we may suppose that the $\mathcal{G}_{j}$ are pairwise distinct and that their ranks are the same.  We write
\begin{equation*}
c_{1}(\mathcal{G}_{j}) = \sum_{i=1}^{k} d_{ij} D_{i}.
\end{equation*}
Since $c_{1}(\mathcal{G}_{j}) \in N_{n-1}(X)_{\mathbb{Z}}$, we know that $d_{ij} \in \frac{1}{m}\mathbb{Z}$.

We first claim that the $d_{ij}$ are bounded above.  Indeed, if not then there would be some index $i$ such that $d_{ij} \to \infty$ as $j \to \infty$.  Then we would have $c_{1}(\mathcal{G}_{j}) \cdot \alpha_{i} \to \infty$ as $j \to \infty$.  But this contradicts the fact that there is an upper bound on the $\alpha_{i}$-slopes of subsheaves of $\mathcal{E}$.

Similarly, the $d_{ij}$ are bounded below.  If not, then there would be some index $i$ such that $d_{ij} \to -\infty$ as $j \to \infty$.  After reindexing and passing to a subsequence, we may assume that $d_{1j}$ is strictly decreasing in $j$.  Among the remaining coefficients, either they are all bounded below as $j$ increases or there is some index which is unbounded.  Rearranging and passing to subsequences inductively, we may suppose that $d_{1j}, d_{2j}, \ldots, d_{tj}$ are strictly decreasing subsequences and the other coefficients $d_{(t+1)j},\ldots,d_{kj}$ are bounded.  Since the latter coefficients are bounded above and below and take values in $\frac{1}{m}\mathbb{Z}$, by passing to a subsequence we may ensure that these other coefficients are constant.  But then since every $\beta_{j}$ is a non-negative linear combination of the $\alpha_{i}$, this would imply that $c_{1}(\mathcal{G}_{j}) \cdot \beta_{j} \leq c_{1}(\mathcal{G}_{1}) \cdot \beta_{j}$ for every $j$.  Since $\mathcal{G}_{1}$ has the same rank as $\mathcal{G}_{j}$, by uniqueness of maximal destabilizing subsheaves we conclude that $\mathcal{G}_{j} = \mathcal{G}_{1}$, a contradiction.

Altogether, we see that the $d_{ij}$ are chosen from a bounded subset of $\frac{1}{m} \mathbb{Z}$.  This means that there are only finitely many possible classes in $N_{n-1}(X)_{\mathbb{Z}}$ coming from a maximal destabilizing subsheaf as we vary over all nef curve classes in $\mathcal{C}$.

Now fix a class $D \in N_{n-1}(X)_{\mathbb{Z}}$ in this finite set and a rank $r$ and consider the subset $\mathcal{C}_{D,r} \subset \mathcal{C}$ of classes $\beta$ such that the corresponding maximal destabilizing subsheaf $\mathcal{G}_{\beta}$ has rank $r$ and has $c_{1}(\mathcal{G}_{\beta}) = D$.  Note that for every $\beta, \beta' \in \mathcal{C}_{D,r}$ the maximal destabilizing subsheaves $\mathcal{G}_{\beta}$, $\mathcal{G}_{\beta'}$ have the same $\beta$-slope and the same $\beta'$-slope.  By the defining property of maximal destabilizing subsheaves, we conclude that $\mathcal{G}_{\beta} \subset \mathcal{G}_{\beta'}$ and conversely.  Thus the subsheaf $\mathcal{G}_{\beta}$ does not change as we vary $\beta \in \mathcal{C}_{D,r}$.   Altogether we see that there is a finite set of $\beta$-maximal destabilizing subsheaves as we vary $\beta \in \mathcal{C}$.
\end{proof}

By combining Theorem \ref{theo:finmanydestabilizing} with an induction on the rank, we obtain:

\begin{corollary} \label{coro:finiteHNfilts}
Let $X$ be a normal projective variety and let $\mathcal{E}$ be a torsion-free sheaf on $X$.  Fix a finite set of nef curve classes $\alpha_{1},\ldots,\alpha_{k} \in \Nef^{n-1}(X)_{\mathbb{Q}}$.  Then there is a finite set of filtrations of $\mathcal{E}$ obtained as $\beta$-Harder-Narasimhan filtrations as we vary $\beta \in \langle \alpha_{1}, \ldots, \alpha_{k} \rangle$.
\end{corollary}

Note that this implies a type of openness for stability regions as in \cite[Section 3]{GKP16}.  Similarly, we have:

\begin{corollary} \label{coro:mumincont}
Let $X$ be a normal projective variety and let $\mathcal{E}$ be a torsion-free sheaf on $X$.  Let $\alpha \in \Nef^{n-1}(X)$ be a curve class in the interior of the cone.  Then $\mu^{min}_{*}(\mathcal{E}): \Nef^{n-1}(X) \to \mathbb{R}$ is continuous at $\alpha$.
\end{corollary}

\begin{proof}
Choose a finite set of rational nef curve classes $\{ \alpha_{i} \}$ whose convex hull $\mathcal{C}$ contains $\alpha$ in its interior.  By Corollary \ref{coro:finiteHNfilts} as we vary $\beta \in \mathcal{C}$ we only get finitely many Harder-Narasimhan filtrations.  For each fixed filtration, the slopes of the quotients depend continuously on $\beta$.  Since $\mu^{min}$ is obtained as a minimum of the slopes of the quotients of the last terms in these finitely many filtrations, we conclude that it is continuous.
\end{proof}

\subsection{Jordan-H\"older filtrations}
An important structural feature of a semistable sheaf is the collection of subsheaves which have the same slope.  In this section we describe some key features of this set.  The nicest properties hold only when $\alpha$ lies in the interior of the nef cone and we will focus on this situation.

\begin{definition}
    Let $X$ be a normal projective variety and let $\alpha \in \Nef^{n-1}(X)$.  Suppose that $\mathcal{E}$ is a semistable sheaf of slope $\mu_{\mathcal{E}}$.  We let $S_{\mathcal{E}}$ denote the set of saturated subsheaves of $\mathcal{E}$ which have slope $\mu_{\mathcal{E}}$.
\end{definition}

By Corollary \ref{coro:maxdestabss} each subsheaf in $S_{\mathcal{E}}$ is semistable.  As we will see, stable sheaves play a particularly important role in analyzing $S_{\mathcal{E}}$.

\begin{lemma} \label{lemm:mapfromstable}
    Let $X$ be a normal projective variety and let $\alpha \in \Nef^{n-1}(X)$.  Suppose that $\mathcal{E}$ is a stable sheaf of slope $\mu_{\mathcal{E}}$.  If $\mathcal{F}$ is a torsion-free sheaf satisfying $\mu^{max}(\mathcal{F}) \leq \mu_{\mathcal{E}}$, then any homomorphism $\phi: \mathcal{E} \to \mathcal{F}$ is either an inclusion or the zero map.
\end{lemma}

\begin{proof}
    Let $\mathcal{K}$ denote the kernel of $\phi$.  Since the image of $\phi$ has slope $\leq \mu^{max}(\mathcal{F})$, Lemma \ref{lemm:exactsequence} shows that either $\mathcal{K}$ is the zero sheaf or $\mu(\mathcal{K}) \geq \mu_{\mathcal{E}}$.  In the latter case, stability implies that $\mathcal{K} \subset \mathcal{E}$ is a full-rank subsheaf.  Thus the image of $\phi$, which is isomorphic to $\mathcal{E}/\mathcal{K}$, is a torsion subsheaf of $\mathcal{F}$, hence $0$.  This shows that $\mathcal{K} = \mathcal{E}$.
\end{proof}

The following definition highlights the most important tool for analyzing $S_{\mathcal{E}}$.

\begin{definition}
    Let $X$ be a normal projective variety and let $\alpha \in \Nef^{n-1}(X)$ be a nef class.  Suppose that $\mathcal{E}$ is a semistable sheaf of slope $\mu_{\mathcal{E}}$.  A Jordan-H\"older filtration for $\mathcal{E}$ is a filtration
    \begin{equation*}
        0 = \mathcal{G}_{0} \subset \mathcal{G}_{1} \subset \ldots \subset \mathcal{G}_{\ell} = \mathcal{E}
    \end{equation*}
    such that each $\mathcal{G}_{i}$ is saturated and each successive quotient $\mathcal{G}_{i}/\mathcal{G}_{i-1}$ is stable of slope $\mu_{\mathcal{E}}$.
\end{definition}

\begin{proposition}[{\cite[Proposition 1.5.2]{HL97}}] \label{prop:jordanholder}
    Let $X$ be a normal projective variety and let $\alpha$ be a class in the interior of $\Nef^{n-1}(X)$.  Suppose that $\mathcal{E}$ is a semistable sheaf of slope $\mu_{\mathcal{E}}$.
    \begin{enumerate}
        \item $\mathcal{E}$ admits a Jordan-H\"older filtration.
        \item In all Jordan-H\"older filtrations of $\mathcal{E}$, the isomorphism type of $\oplus_{i=1}^{\ell} \widehat{\mathcal{G}_{i}/\mathcal{G}_{i-1}}$ is always the same.
    \end{enumerate}
\end{proposition}

\begin{proof}
    (1) Let $\mathcal{G}_{1}$ be any non-zero subsheaf in $S_{\mathcal{E}}$ with the minimum possible rank.  Then $\mathcal{G}_{1}$ must be stable.  If $\mathcal{E} = \mathcal{G}_{1}$, we are done.  Otherwise $\mathcal{E}/\mathcal{G}_{1}$ still has slope $\mu_{\mathcal{E}}$ and thus is still semistable.  We conclude by induction on the rank.

    (2) The proof is by induction on the rank $r$ of $\mathcal{E}$.  The base case $r=1$ is trivial.  For the induction step, suppose that $\{\mathcal{G}_{i}\}$ and $\{ \mathcal{G}_{j}' \}$ are Jordan-H\"older filtrations of $\mathcal{E}$.  Let $k$ be the minimal index such that $\mathcal{G}_{1}' \subset \mathcal{G}_{k}$ as subsheaves of $\mathcal{E}$.  By Lemma \ref{lemm:mapfromstable}, the induced map $\psi: \mathcal{G}_{1}' \to \mathcal{G}_{k}/\mathcal{G}_{k-1}$ is either the zero morphism or an inclusion, and the former possibility would violate the minimality of $k$.  Since $\mathcal{G}_{k}/\mathcal{G}_{k-1}$ is stable, the equality of slopes implies that $\psi$ is an inclusion of sheaves of equal rank.  Since $\alpha$ lies in the interior of the nef cone, Proposition \ref{prop:c1properties}.(3) shows that $\psi$ is an isomorphism in codimension $1$.  Thus there is a big open subset $U$ such that $\mathcal{G}_{k}|_{U} \cong \mathcal{G}_{1}'|_{U} \oplus \mathcal{G}_{k-1}|_{U}$ as subsheaves of $\mathcal{E}$.  In particular $\widehat{\mathcal{G}_{1}'} \cong \widehat{\mathcal{G}_{k}/\mathcal{G}_{k-1}}$.
    
    The sheaf $\mathcal{E}/\mathcal{G}_{1}'$ has two Jordan-H\"older filtrations.  The first is the filtration $\{ \mathcal{G}_{j}'/\mathcal{G}_{1}'\}$.  We construct the second from $\{ \mathcal{G}_{i} \}$ as follows.  For $i < k$, we set $\mathcal{H}_{i}$ to be the saturation in $\mathcal{E}/\mathcal{G}'_{1}$ of $(\mathcal{G}'_{1} + \mathcal{G}_{i})/\mathcal{G}'_{1}$.  For $i > k$, we set $\mathcal{H}_{i} = \mathcal{G}_{i}/\mathcal{G}'_{1}$.  Note that for every $i \neq k,k+1$ we have $\mathcal{G}_{i}/\mathcal{G}_{i-1}|_{U} \cong \mathcal{H}_{i}/\mathcal{H}_{i-1}|_{U}$ after possibly replacing $U$ by a big open subset.  In particular, all the quotients $\mathcal{H}_{i}/\mathcal{H}_{i-1}$ are stable and their reflexive hulls agree with that of $\mathcal{G}_{i}/\mathcal{G}_{i-1}$.  We have a similar relationship between $\mathcal{G}_{k+1}/\mathcal{G}_{k}$ and $\mathcal{H}_{k+1}/\mathcal{H}_{k-1}$.
    
    By applying the induction assumption to $\mathcal{E}/\mathcal{G}_{1}'$, we obtain the desired conclusion for $\mathcal{E}$.
\end{proof}

Note that there may be many different Jordan-H\"older filtrations and that $S_{\mathcal{E}}$ need not be finite.
There are several factors determining how $S_{\mathcal{E}}$ is constructed from the Jordan-H\"older factors of $\mathcal{E}$.  First, two stable factors could have a non-trivial extension group.  Second, for a given extension and a given subsheaf there may be many ways to embed the subsheaf into the extension.  Despite these complications, we can give some structure to $S_{\mathcal{E}}$.

\begin{lemma} \label{lemm:sameslopeint}
    Let $X$ be a normal projective variety and let $\alpha \in \Nef^{n-1}(X)$.  Suppose that $\mathcal{E}$ is a semistable sheaf of slope $\mu_{\mathcal{E}}$.  For any two sheaves $\mathcal{F}_{1},\mathcal{F}_{2} \in S_{\mathcal{E}}$:
    \begin{enumerate}
        \item $\mathcal{F}_{1} \cap \mathcal{F}_{2}$ is either in $S_{\mathcal{E}}$ or is the zero sheaf.
        \item The saturation of $\mathcal{F}_{1} + \mathcal{F}_{2}$ is in $S_{\mathcal{E}}$.
    \end{enumerate}
\end{lemma}

\begin{proof}
    Since we have an inclusion $\mathcal{E}/(\mathcal{F}_{1} \cap \mathcal{F}_{2}) \hookrightarrow (\mathcal{E}/\mathcal{F}_{1}) \oplus (\mathcal{E}/\mathcal{F}_{2})$ we see that $\mathcal{F}_{1} \cap \mathcal{F}_{2}$ is a saturated subsheaf.  From the exact sequence
    \begin{equation*}
        0 \to \mathcal{F}_{1} \cap \mathcal{F}_{2} \to \mathcal{F}_{1} \oplus \mathcal{F}_{2} \to \mathcal{F}_{1} + \mathcal{F}_{2} \to 0
    \end{equation*}
    we see that $\mu_{\mathcal{E}}$ is a weighted average of $\mu(\mathcal{F}_{1} \cap \mathcal{F}_{2})$ and $\mu(\mathcal{F}_{1} + \mathcal{F}_{2})$.  If $\mathcal{F}_{1} \cap \mathcal{F}_{2}$ is the zero sheaf, we are done.  Otherwise, both outer terms are subsheaves of $\mathcal{E}$ and so have slope at most $\mu_{\mathcal{E}}$.  We conclude that both sheaves must have slope exactly $\mu_{\mathcal{E}}$.  This implies that the saturation of $\mathcal{F}_{1} + \mathcal{F}_{2}$ also has slope $\mu_{\mathcal{E}}$.
\end{proof}

\begin{remark}
    The proof of Lemma \ref{lemm:sameslopeint} shows that $\mathcal{F}_{1} + \mathcal{F}_{2}$ has the same slope as its saturation.  In particular, if $\alpha$ lies in the interior of $\Nef^{n-1}(X)$ then the torsion subsheaf of $\mathcal{E}/(\mathcal{F}_{1}+\mathcal{F}_{2})$ is supported in codimension $\geq 2$.
\end{remark}

The next definition identifies another useful tool for describing the relationship between stability and semistability.

\begin{definition}
    Let $X$ be a normal projective variety and let $\alpha \in \Nef^{n-1}(X)$.  We say that a sheaf $\mathcal{E}$ is polystable if it is a direct sum of stable sheaves which all have the same slope.
\end{definition}

Just as with Jordan-H\"older filtrations, the theory of polystable sheaves requires some extra work to account for codimension $\geq 2$ subsets.  A convenient way to handle these issues is to work with reflexive sheaves as in the following proposition.

\begin{proposition}[{\cite[Proposition 1.5.5]{HL97}}] 
    Let $X$ be a normal projective variety and let $\alpha$ be a class in the interior of $\Nef^{n-1}(X)$.  Suppose that $\mathcal{E}$ is a reflexive semistable sheaf of slope $\mu_{\mathcal{E}}$.  Amongst all polystable subsheaves of $\mathcal{E}$ with slope $\mu_{\mathcal{E}}$, there is a unique maximal element $\mathcal{P}$.
\end{proposition}

The sheaf $\mathcal{P}$ is called the socle of $\mathcal{E}$.

\begin{proof}
    We inductively construct polystable subsheaves $\mathcal{F}_{j}$ in $S_{\mathcal{E}}$ as follows.  Proposition \ref{prop:jordanholder}.(1) shows that $\mathcal{E}$ admits a saturated stable subsheaf $\mathcal{F}_{1}$ of slope $\mu_{\mathcal{E}}$.  Inductively, if every stable subsheaf of slope $\mu_{\mathcal{E}}$ is contained in $\mathcal{F}_{j}$, we set $\mathcal{P} = \mathcal{F}_{j}$.  Otherwise, let $\mathcal{F}_{j+1}$ be the saturation of the sum of $\mathcal{F}_{j}$ and a stable saturated sheaf $\mathcal{G}_{j}$ of slope $\mu_{\mathcal{E}}$ that is not contained in $\mathcal{F}_{j}$.  By Lemma \ref{lemm:sameslopeint}.(2) $\mathcal{F}_{j+1}$ is still in $S_{\mathcal{E}}$.  Furthermore, by Lemma \ref{lemm:sameslopeint}.(1) we have $\mathcal{F}_{j} \cap \mathcal{G}_{j}$ is a proper saturated subsheaf of $\mathcal{G}_{j}$ with the same slope, hence $0$.  We conclude that we have an inclusion $\mathcal{F}_{j} \oplus \mathcal{G}_{j} \to \mathcal{F}_{j+1}$ that is an isomorphism away from a closed subset of codimension $\geq 2$.  Since $\mathcal{F}_{j},\mathcal{G}_{j},\mathcal{F}_{j+1}$ are all saturated in $\mathcal{E}$, hence reflexive, this map is an isomorphism.  We deduce that $\mathcal{F}_{j+1}$ is still polystable.
    
    Since the rank of a subsheaf of $\mathcal{E}$ is bounded above, this process must terminate with a polystable subsheaf $\mathcal{P}$ that contains every stable subsheaf of slope $\mu_{\mathcal{E}}$.  Thus $\mathcal{P}$ satisfies the desired maximality property.
\end{proof}

\section{Semistability under geometric operations} \label{sect:geomop}

We next turn to the behavior of semistability of torsion-free sheaves under geometric operations: tensor products and finite pullbacks.

For smooth projective varieties in characteristic $0$, the behavior of semistability under these operations is well-understood.  The situation for normal projective varieties is no different.

\begin{theorem} \label{theo:char0props}
    Let $X$ be a normal projective variety of dimension $n$ over an algebraically closed field of characteristic $0$.  Let $\alpha \in \Nef^{n-1}(X)$.
    \begin{enumerate}
        \item  Suppose that $\mathcal{E}_{1},\mathcal{E}_{2}$ are semistable torsion-free sheaves.  Then $\mathcal{E}_{1} \widehat{\otimes} \mathcal{E}_{2}$ is semistable.
        \item Let $f: Y \to X$ be a generically finite surjective morphism from a normal projective variety $Y$.  For every torsion-free coherent sheaf $\mathcal{E}$ on $X$, $\mathcal{E}$ is $\alpha$-semistable if and only if $(f^{*}\mathcal{E})_{tf}$ is $f^{*}\alpha$-semistable.
    \end{enumerate}
\end{theorem}

\begin{proof}
    (1) When $X$ is a smooth projective variety, this is proved in \cite[Theorem 5.1]{CP11} or \cite[Theorem 4.2]{GKP16} (over $\mathbb{C}$, and thus by a standard argument, over an arbitrary algebraically closed field of characteristic $0$).  The general case follows from Lemma \ref{lemm:slopebirational}.
    
    (2) We first show the reverse implication.  Suppose $\mathcal{E}$ is not semistable and let $\mathcal{F}$ be a maximal destabilizing subsheaf.  Then $(f^{*}\mathcal{F})_{tf} \subset (f^{*}\mathcal{E})_{tf}$ is a destabilizing subsheaf; indeed, the slopes of both sheaves are simply multiplied by $\deg(f)$.

    Conversely, we first suppose that $f: Y \to X$ has the property that $K(Y)/K(X)$ is Galois with Galois group $G$.  Let $Y'$ denote the normalization of $X$ in $K(Y)$, so that $f$ factors as a birational $g: Y \to Y'$ followed by a finite $h: Y' \to X$.  By Lemma \ref{lemm:slopebirational} it suffices to prove the desired statement for $h$.

    Suppose for a contradiction that $(h^{*}\mathcal{E})_{tf}$ is not semistable.  Note that $G$ acts on $h^{*}\mathcal{E}$, and thus on its torsion-free quotient. The maximal destabilizing subsheaf $\mathcal{F} \subset (h^{*}\mathcal{E})_{tf}$  is $G$-invariant.  Let $V \subset Y'$ be the big open subset where $h$ is flat.  By applying \cite[Lemma 2.4]{Weissmann24} along $V$, we obtain a saturated subsheaf $\mathcal{G} \subset \mathcal{E}$ such that $h^{*}(\mathcal{G}|_{h(V)}) = \mathcal{F}|_{V}$.  Since
    \begin{equation*}
        \mu_{h^{*}\alpha}((h^{*}\mathcal{E})_{tf}) < \mu_{h^{*}\alpha}(\mathcal{F}) = \mu_{h^{*}\alpha}((h^{*}\mathcal{G})_{tf}) = \deg(h) \mu_{\alpha}(\mathcal{G}) \leq \deg(h) \mu_{\alpha}(\mathcal{E}) = \mu_{h^{*}\alpha}((h^{*}\mathcal{E})_{tf})
    \end{equation*}
    we obtain a contradiction.

    In the general case, as before we can take a Stein factorization and reduce the problem to the finite part using Lemma \ref{lemm:slopebirational}.  Then by the first paragraph of the proof of (2) it suffices to prove semistability on any higher finite cover.  In particular we may pass to a Galois cover, where the result has already been proved.
\end{proof}

In characteristic $p$, we can no longer expect such nice properties.  It turns out that the theory for arbitrary nef curve classes is essentially no different than the theory for complete intersection classes as developed in \cite{Langer15}.  In the following subsections we review the foundations of this theory.

\subsection{Asymptotic Frobenius pullbacks}

From now on we work over an algebraically closed field $\mathbb{K}$ of positive characteristic.  We first address the case of separable maps.

\begin{theorem} \label{theo:separablepullbackisss}
    Let $X$ be a normal projective variety of dimension $n$ over an algebraically closed field of characteristic $p$.  Let $\alpha \in \Nef^{n-1}(X)$.
    
    Let $f: Y \to X$ be a generically finite surjective separable morphism from a normal projective variety $Y$.  For every torsion-free coherent sheaf $\mathcal{E}$ on $X$, $\mathcal{E}$ is $\alpha$-semistable if and only if $(f^{*}\mathcal{E})_{tf}$ is $f^{*}\alpha$-semistable.
\end{theorem}

\begin{proof}
    The proof is the same as the proof of Theorem \ref{theo:char0props}.(2).
\end{proof}

We next turn to purely inseparable maps.  In general it is not true that Frobenius pullback preserves semistability of torsion-free sheaves.  The following definition highlights the torsion-free sheaves with this special property.

\begin{definition}
    Let $X$ be a normal projective variety over an algebraically closed field of characteristic $p$ and let $\alpha \in \Nef^{n-1}(X)$ be a nef curve class.  A non-zero torsion-free sheaf $\mathcal{E}$ is strongly $\alpha$-semistable (respectively strongly $\alpha$-stable) if $(F^{e*}_{abs}\mathcal{E})_{tf}$ is $\alpha$-semistable (respectively $\alpha$-stable) for every positive integer $e$.  We will drop the prefix ``$\alpha$'' when it is clear from the context.
\end{definition}

By Remark \ref{rema:bigopen} we may equally well require that $F^{e\widehat{*}}_{abs}\mathcal{E}$ is $\alpha$-semistable for every $e \geq 1$, and we will usually use this formulation instead.  

\begin{lemma}
    Let $X$ be a normal projective variety over an algebraically closed field of characteristic $p$ and let $\alpha \in \Nef^{n-1}(X)$ be a nef curve class.  Suppose that $\mathcal{E}$ is a strongly $\alpha$-semistable torsion-free sheaf on $X$.  Then for any finite purely inseparable morphism $f: Y \to X$ of normal projective varieties, we have that $f^{\widehat{*}}\mathcal{E}$ is $f^{*}\alpha$-semistable.
\end{lemma}

The same statement holds for strong stability.

\begin{proof}
      Choose $e$ large enough so that $F^{e}_{rel}$ factors as $f \circ g$ for some finite morphism $g: X^{(-e)} \to Y$.  Suppose for a contradiction that there were a destabilizing subsheaf $\mathcal{G} \subset f^{\widehat{*}}\mathcal{E}$.  We would then obtain a destabilizing subsheaf $g^{\widehat{*}}\mathcal{G} \subset F^{e\widehat{*}}_{rel}\mathcal{E}$.  However, the abstract isomorphism $i: X \xrightarrow{\cong} X^{(-e)}$ identifies the sheaf $F^{e\widehat{*}}_{abs}\mathcal{E}$ with the sheaf $F^{e\widehat{*}}_{rel}\mathcal{E}$ and identifies $\alpha$ and $F_{rel}^{e*}\alpha$ up to a positive scaling factor.  Since $i$ is $\mathbb{K}$-semilinear, it identifies the degree function on $0$-cycles for $X$ and $X^{(-e)}$ and thus preserves intersection products.  Thus we see that $i^{*}(g^{\widehat{*}}\mathcal{G})$ destabilizes $F^{e\widehat{*}}_{abs}\mathcal{E}$, a contradiction to strong semistability.
\end{proof}

\cite[Theorem 2.7]{Langer04} shows that the Harder-Narasimhan filtrations of iterated Frobenius pullbacks of $\mathcal{E}$ eventually stabilize.  For the proof it is convenient to use the Harder-Narasimhan polygon of $\mathcal{E}$.

\begin{definition}
    Let $X$ be a normal projective variety over an algebraically closed field of characteristic $p$ and let $\alpha \in \Nef^{n-1}(X)$ be a nef curve class.  Consider a non-zero torsion-free sheaf $\mathcal{E}$ of rank $r$ with Harder-Narasimhan filtration
    \begin{equation*}
        0 = \mathcal{F}_{0} \subset \ldots \subset \mathcal{F}_{s} = \mathcal{E}.
    \end{equation*}
    The Harder-Narasimhan polygon $\HNP(\mathcal{E})$ is the union of the line segments connecting the points $(\rk(\mathcal{F}_{i-1}),c_{1}(\mathcal{F}_{i-1}) \cdot \alpha)$ and $(\rk(\mathcal{F}_{i}),c_{1}(\mathcal{F}_{i}) \cdot \alpha)$ for $1 \leq i \leq s$.  We can equivalently think of this polygon as the graph of a piecewise linear function $f: [0,r] \to \mathbb{R}$.  We say that the $x$-value $i \in \{1,\ldots,r-1\}$ is a vertex if $f$ fails to be linear at $i$.

    For any positive integer $e$, we define the $e$th Harder-Narasimhan polygon $\HNP_{e}(\mathcal{E})$ to be the Harder-Narasimhan polygon $\HNP(F^{e\widehat{*}}_{abs}\mathcal{E})$ rescaled by $1/p^{e}$ in the second coordinate.  We let $f_{e}: [0,r] \to \mathbb{R}$ denote the corresponding piecewise linear function.
\end{definition}

Note that $f_{e}(x) \leq f_{e+1}(x)$ for every $x \in [0,r]$. The following lemma allows us to take a limit as $e \to \infty$.

\begin{lemma}
    Let $X$ be a normal projective variety over an algebraically closed field of characteristic $p$ and let $\alpha \in \Nef^{n-1}(X)$ be a nef curve class.  Let $\mathcal{E}$ be a non-zero torsion-free sheaf.  Then there is a constant $C$ such that for every positive integer $e$ the polygon $\HNP_{e}(\mathcal{E})$ lies below the line $y=C$.
\end{lemma}

\begin{proof}
    Choose an injection $\mathcal{E} \to \mathcal{O}_{X}(H)^{\oplus k}$ where $H$ is a sufficiently positive Cartier divisor.  Then we also have injections $F^{e\widehat{*}}_{abs}\mathcal{E} \to \mathcal{O}_{X}(p^{e}H)^{\oplus k}$.  The latter bundle is $\alpha$-semistable by Lemma \ref{lemm:semistabpreserved}.  The desired statement follows by setting $C = (\rk \mathcal{E}) H \cdot \alpha$.
\end{proof}

\begin{definition}
    Let $X$ be a normal projective variety over an algebraically closed field of characteristic $p$ and let $\alpha \in \Nef^{n-1}(X)$ be a nef curve class.  Consider a non-zero torsion-free sheaf $\mathcal{E}$.  We define the Harder-Narasimhan polygon $\HNP_{\infty}(\mathcal{E})$ as the graph of the piecewise linear function $f_{\infty} = \lim_{e \to \infty} f_{e}$.
\end{definition}

\begin{theorem}[{\cite[Theorem 2.7]{Langer04}}] \label{theo:asymsss}
    Let $X$ be a normal projective variety over an algebraically closed field of characteristic $p$ and let $\alpha \in \Nef^{n-1}(X)$ be a nef curve class.

    Let $\mathcal{E}$ be a non-zero torsion-free sheaf on $X$.  There is some positive integer $e_{0}$ such that for every $e \geq e_{0}$ the $\alpha$-Harder-Narasimhan filtration of $F_{abs}^{e\widehat{*}}\mathcal{E}$ is the $F_{abs}^{(e-e_{0})}$-reflexive pullback of the $\alpha$-Harder-Narasimhan filtration of $F_{abs}^{e_{0}\widehat{*}}\mathcal{E}$.
    
    In particular, every successive quotient in the Harder-Narasimhan filtration of $F_{abs}^{e_{0}\widehat{*}}\mathcal{E}$ is strongly semistable.
\end{theorem}

\begin{proof}
    The proof is by induction on the rank $r$ of $\mathcal{E}$.  The base case when the rank is $1$ is obvious.  For the induction step, it suffices to consider the case when $\mathcal{E}$ is not strongly semistable.  In other words, the set of vertices for $f_{\infty}$ is non-empty.  Note that $f_{\infty}$ is a piecewise linear function whose slopes are decreasing as $x$ increases.  From now on we let $b$ denote the $x$-value of the first vertex and let $\mu_{diff}$ denote the difference in slopes between the first two linear pieces of $f_{\infty}$.  
    
    Consider $f_{e}$ for some $e \geq 1$.  If $K_{e} \subset \{1,\ldots,r-1\}$ denotes the set of vertices, then the Harder-Narasimhan filtration has the form
    \begin{equation*}
        0 = \mathcal{F}_{e,0} \subset \mathcal{F}_{e,k_{1}} \subset \ldots \subset \mathcal{F}_{e,k_{t_{e}}} \subset \mathcal{F}_{e,r} = F^{e \widehat{*}}_{abs}\mathcal{E}
    \end{equation*}
    where $k_{1},\ldots,k_{t_{e}}$ are the elements of $K_{e}$.
    
    Note that as $e$ increases, $f_{e}$ converges monotonically and pointwise to $f_{\infty}$.  Thus for any fixed $\epsilon$, we have $\sup_{x \in [0,r]} f_{\infty}(x) - f_{e}(x) < \epsilon$ for all sufficiently large $e$.  Note also that there must be infinitely many values of $e$ such that $f_{e}$ has a vertex at $b$, each corresponding to a term $\mathcal{F}_{e,b}$ in the Harder-Narasimhan filtration of $F^{e \widehat{*}}_{abs} \mathcal{E}$.
    
    First suppose that for some $e'$, $f_{e'}$ has a vertex at $b$ and the first line segment of the rescaled polygon $\frac{1}{p^{e'}}\HNP_{\infty}(\mathcal{F}_{e',b})$ lies on the first line segment of $\HNP_{\infty}(\mathcal{E})$.  By the induction assumption, there is a positive integer $q$ such that $\HNP_{\infty}(\mathcal{F}_{e',b}) = \HNP_{q}(\mathcal{F}_{e',b})$.  Let $\mathcal{G}$ be the maximal destabilizing subsheaf of $F^{q\widehat{*}}_{abs}\mathcal{F}_{e',b}$.  Our vertex assumption implies that $\mathcal{G}$ is strongly semistable and that for every $e \geq q+e'$ the sheaf $F^{(e-q-e')\widehat{*}}_{abs}\mathcal{G}$ is a subsheaf of $F^{e\widehat{*}}_{abs}\mathcal{E}$ with the maximal possible slope. Then we can apply the induction assumption to the quotient $F^{(q+e')\widehat{*}}_{abs}\mathcal{E}/\mathcal{G}$ to deduce the desired property.

    Next suppose for a contradiction that for every $e'$ such that $f_{e'}$ has a vertex at $b$ the first line segment of $\frac{1}{p^{e'}}\HNP_{\infty}(\mathcal{F}_{e',b})$ lies strictly below the first line segment of $\HNP_{\infty}(\mathcal{E})$.  Define
    \begin{equation*}
        \epsilon = \frac{\mu_{diff}}{2r}
    \end{equation*}
    and suppose that $e$ lies in the range where $\sup_{x \in [0,r]} f_{\infty}(x) - f_{e}(x) < \epsilon$ and that $f_{e}$ has a vertex at $b$.  By our contradiction assumption, there is a $k_{j} \in K_{e} \cup \{ r \}$ with $k_{j} > b$ and an $\ell > e$ such that
    \begin{equation*}
        \mu^{max}_{\alpha}(F^{(\ell-e)\widehat{*}}_{abs} (\mathcal{F}_{e,k_{j}}/ \mathcal{F}_{e,k_{j-1}})) > \mu^{max}_{\alpha}(F^{(\ell-e)\widehat{*}}_{abs} \mathcal{F}_{e,b}).
    \end{equation*}
    There is a map $F^{(\ell-e)\widehat{*}}_{abs} \mathcal{F}_{e,k_{j}}/ F^{(\ell-e)\widehat{*}}_{abs}\mathcal{F}_{e,k_{j-1}} \to F^{(\ell-e)\widehat{*}}_{abs} (\mathcal{F}_{e,k_{j}}/ \mathcal{F}_{e,k_{j-1}})$ that is an isomorphism away from a codimension $\geq 2$ set so that the Harder-Narasimhan filtrations of the two sheaves can naturally be identified. Let $\mathcal{G} \subset F^{(\ell-e)\widehat{*}}_{abs} \mathcal{F}_{e,k_{j}}$ denote the preimage of the maximal destabilizing subsheaf of $F^{(\ell-e)\widehat{*}}_{abs} \mathcal{F}_{e,k_{j}}/ F^{(\ell-e)\widehat{*}}_{abs}\mathcal{F}_{e,k_{j-1}}$.
    Then the graph of $f_{\infty}$ lies above both the graph of $f_{e}$ and the point $(\rk(\mathcal{G}),c_{1}(\mathcal{G}) \cdot \alpha/p^{\ell})$.  Let $g: [0,r] \to \mathbb{R}$ denote the minimal piecewise linear function with these two properties that only fails to be linear at integer values of $x$.  We have
    \begin{align*}
        \int_{0}^{r} f_{\infty}(x) - f_{e}(x) \, dx & \geq \int_{0}^{r} g(x) - f_{e}(x) \, dx \\
        & \geq \frac{1}{2} \left(\frac{\mu_{\alpha}(\mathcal{G}/F^{(\ell-e)\widehat{*}}_{abs}\mathcal{F}_{e,k_{j-1}})}{p^{\ell}} - \frac{\mu
        _{\alpha}(F^{(\ell-e)\widehat{*}}_{abs}(\mathcal{F}_{e,k_{j}}/\mathcal{F}_{e,k_{j-1}}))}{p^{\ell}}  \right) \\
        & > \frac{1}{2} \left( \frac{\mu_{\alpha}(\mathcal{F}_{e,b})}{p^{e}} - \frac{\mu_{\alpha}(\mathcal{F}_{e,k_{j}}/\mathcal{F}_{e,k_{j-1}})}{p^{e}} \right) \\
        & \geq \frac{1}{2} (f_{\infty}(1) - (f_{\infty}(k_{j}) - f_{\infty}(k_{j}-1)) - 2 \epsilon) \\
        & \geq \frac{\mu_{diff}}{2} - \epsilon.
    \end{align*}
    On the other hand, by construction we have
    \begin{align*}
        \int_{0}^{r} f_{\infty}(x) - f_{e}(x) \, dx & = \left( \int_{0}^{1} + \int_{1}^{r-1} + \int_{r-1}^{r} \right) f_{\infty}(x) - f_{e}(x) \, dx 
        \\ & \leq \frac{\epsilon}{2} + (r-2) \epsilon + \frac{\epsilon}{2} \\
        & = (r-1)\epsilon.
    \end{align*}
    Combined these yield a contradiction, finishing the proof.
\end{proof}

It will be useful to isolate some of the information encapsulated in $\HNP_{\infty}(\mathcal{E})$.  Note that we have
\begin{equation*}
    \mu_{\alpha}^{min}(F^{\widehat{*}}_{abs}\mathcal{E})/p \leq \mu_{\alpha}^{min}(\mathcal{E}) \leq \mu_{\alpha}^{max}(\mathcal{E}) \leq \mu_{\alpha}^{max}(F^{\widehat{*}}_{abs}\mathcal{E})/p.
\end{equation*}
We define
\begin{equation*}
    L^{max}(\mathcal{E}) = \limsup_{k \to \infty} \frac{\mu^{max}_{\alpha}(F^{k\widehat{*}}_{abs}\mathcal{E})}{p^{k}} \qquad \qquad L^{min}(\mathcal{E}) = \liminf_{k \to \infty} \frac{\mu^{min}_{\alpha}(F^{k\widehat{*}}_{abs}\mathcal{E})}{p^{k}} 
\end{equation*}
Note that the sequences defining the invariants are monotonic and even eventually constant by Theorem \ref{theo:asymsss}.  For every $e \geq 0$ we have
\begin{equation*}
    L^{min}(\mathcal{E}) \leq \mu_{\alpha}^{min}(F^{e\widehat{*}}_{abs}\mathcal{E})/p^{e}  \qquad \textrm{and} \qquad \mu_{\alpha}^{max}(F^{e\widehat{*}}_{abs}\mathcal{E})/p^{e} \leq L^{max}(\mathcal{E}).
\end{equation*}

\subsection{Tensor products in characteristic p} \label{sect:tensorprods}

Let $X$ be a normal projective variety over an algebraically closed field of characteristic $p$.  Suppose $\mathcal{E}$ is a non-zero torsion-free sheaf on $X$.  Let $U \subset X$ be a big open subset such that $\mathcal{E}|_{U}$ is locally free of rank $r$.  Then the frame bundle $\mathcal{G}$ associated to $\mathcal{E}|_{U}$ is a principal $GL(r)$-bundle.  The theory of slope semistability of $\mathcal{E}$ is closely related to the theory of slope semistability of $\mathcal{G}$,  which we now quickly recall.  This theory is developed in \cite{RR84}, \cite{IMP03}, \cite{Langer04b}, \cite{Langer05}, and \cite{CH06}.

\begin{notation}
    Suppose that $G$ is an algebraic group, $X$ is a variety with a right $G$-action, and $Y$ is a variety with a left $G$-action.  We define $X \times^{G} Y = (X \times Y)/G$ where $G$ acts via the right action $(x,y) \cdot g = (xg,g^{-1}y)$.

      A parabolic subgroup will always implicitly be assumed to be reduced.
\end{notation}

\subsubsection{Roots and heights}

Let $G$ be a connected reductive group over $\mathbb{K}$.  The radical $Z_{0}(G)$ of $G$ coincides with the reduced connected component of the center of $G$ that contains the identity.

Suppose we fix a choice $T \subset B$ of a maximal torus $T$ and a Borel subgroup $B$ in $G$.  Recall that $B$ corresponds to a choice of positive roots $\Phi(G,T)^{+}$ inside the set of roots $\Phi(G,T) \subset \mathrm{Hom}(T,\mathbb{G}_{m})$, where $\mathrm{Hom}(T,\mathbb{G}_{m})$ is the lattice of group homomorphisms.  Letting $\Delta \subset \Phi(G,T)^{+}$ denote the set of simple roots, there is a bijection between subsets $I \subset \Delta$ and parabolic subgroups $P_{I}$ containing $B$.

Suppose $\chi: P \to \mathbb{G}_{m}$ is a character of a parabolic subgroup $P \supset B$.  Then $\chi$ is determined by its restriction to $T$.  The character $\chi$ is said to be dominant if it vanishes on $Z_{0}(G)$ and its pairing against the coroot $\alpha^{\vee}$ is non-negative for every $\alpha \in \Delta$ (or equivalently, every $\alpha \in \Delta \backslash I$ where $I$ corresponds to $P$).  The notion of dominance is independent of the choice of both $T$ and $B$ contained in $P$.

Each character $\chi$ of $P$ produces a $G$-equivariant line bundle $L_{\chi}$ on the projective variety $G/P$ via the prescription $L_{\chi} = G \times^{P} \mathbb{A}^{1}_{\chi}$.  When $\chi$ is dominant, the line bundle $L^{-1}_{\chi}$ is globally generated.

We define the height of a character $\chi: T \to \mathbb{G}_{m}$ by the formula
\begin{equation*}
    \mathrm{ht}(\chi) = \frac{1}{2} \sum_{\alpha \in \Phi(G,T)^{+}} \langle \chi,\alpha^{\vee} \rangle.
\end{equation*}
Note that if $\chi$ is dominant then its height is non-negative.

\subsubsection{Principal $G$-bundles}  
All principal $G$-bundles will admit a right action and will be defined with respect to the fppf topology.

\begin{definition}
    Let $Y$ be a smooth variety. Let $G$ be a connected reductive group over $\mathbb{K}$ and suppose that $\mathcal{G}$ is a principal $G$-bundle over $Y$.  Given a subgroup $H \subset G$, we say that a principal $H$-bundle $\mathcal{H}$ over $Y$ is a reduction of $\mathcal{G}$ from $G$ to $H$ if there is an isomorphism $\mathcal{G} \cong \mathcal{H} \times^{H} G$.

    There is a bijection between reductions of $\mathcal{G}$ from $G$ to $H$ and sections $\sigma: Y \to \mathcal{G}/H$. This bijection identifies a section $\sigma$ with the reduction obtained by pulling back the principal $H$-bundle $\mathcal{G} \to \mathcal{G}/H$ by $\sigma$.
\end{definition}

For a principal $G$-bundle $\mathcal{G}$ and a character $\chi$ of a parabolic subgroup $P$, we define the line bundle $L_{\chi}$ on $\mathcal{G}/P$ as $L_{\chi} = \mathcal{G} \times^{P} \mathbb{A}^{1}_{\chi}$. 

\begin{definition}
    Let $X$ be a normal projective variety over an algebraically closed field of characteristic $p$ and let $\alpha \in \Nef^{n-1}(X)$ be a nef curve class.  Let $G$ be a connected reductive group and suppose that $\mathcal{G}$ is a principal $G$-bundle over a big open set $U \subset X$ contained in the smooth locus of $X$.  We say that $\mathcal{G}$ is $\alpha$-semistable (respectively $\alpha$-stable) if for every big open subset $U' \subset U$, for every reduction $\sigma: U' \to \mathcal{G}/P$ from $G$ to a parabolic subgroup $P$, and for every non-trivial dominant character $\chi$ on $P$, the reflexive extension $\widetilde{\sigma^{*}L_{\chi}}$ of $\sigma^{*}L_{\chi}$ satisfies
    \begin{equation} \label{eq:gss}
        c_{1}(\widetilde{\sigma^{*}L_{\chi}}) \cdot \alpha \leq 0 \qquad (\textrm{respectively }<0).
    \end{equation}
    We say that $\mathcal{G}$ is strongly $\alpha$-semistable (respectively strongly $\alpha$-stable) if $F^{e*}_{abs}\mathcal{G}|_{U}$ is $\alpha$-semistable (respectively $\alpha$-stable) for every $e \geq 1$.  We will frequently drop the prefix ``$\alpha$'' when it is clear from context.
\end{definition}

The following lemma connects the notions of stability in the settings of coherent sheaves and principal $G$-bundles.

\begin{lemma}[{\cite[Lemma A.2.2]{Langer04b}}] \label{lemm:equivalentss}
    Let $X$ be a normal projective variety over an algebraically closed field of characteristic $p$ and let $\alpha \in \Nef^{n-1}(X)$ be a nef curve class.  Let $\mathcal{E}$ be a non-zero torsion-free sheaf of rank $r$ on $X$ and let $U \subset X$ be a big open set contained in the smooth locus of $X$ such that $\mathcal{E}|_{U}$ is locally free.
    
    Let $\mathcal{G} = \underline{\mathrm{Isom}}_{U}(\mathcal{O}_{U}^{\oplus r},\mathcal{E})$ be the principal $GL(r)$-bundle on $U$ obtained by taking the frame bundle of $\mathcal{E}|_{U}$.  Then $\mathcal{E}$ is semistable if and only if $\mathcal{G}$ is semistable.
\end{lemma}

The analogous statement is true for stability, strong semistability, and strong stability by the same proof.

\begin{proof}
    We first show that semistability of $\mathcal{G}$ implies semistability of $\mathcal{E}$.  Suppose we are given a non-zero saturated torsion-free subsheaf $\mathcal{F} \subset \mathcal{E}$ of rank $s$.  On some big open subset $U' \subset U$, the restriction $\mathcal{F}|_{U'} \subset \mathcal{E}|_{U'}$ is a saturated inclusion of locally free sheaves.  Fix a rank $s$ subspace $S \subset \mathbb{K}^{r}$ and let $P \subset GL(r)$ be the maximal parabolic subgroup that stabilizes $S$.  Then we can define a principal $P$-bundle $\mathcal{P}$ on $U'$ via the functor-of-points description $T \mapsto \{ \phi \in \mathcal{G}(T) | \phi(S_{T}) = \mathcal{F}_{T} \}$.  Note that $\mathcal{P}$ is a reduction of $\mathcal{G}$ from $G$ to $P$ and thus defines a section $\sigma: U' \to \mathcal{G}/P$.
    
    Any character $\chi$ on $P$ is determined by its restriction to $GL(s) \times GL(r-s)$ and thus will be determined by its restrictions to the two factors.  Since $P$ is maximal, up to scaling there is a unique character $\chi$ that vanishes along $Z_{0}(GL(r))$.  The dominant character $\chi = (r-s) \det_{s} - s \det_{r-s}$ is a natural choice.  Then $c_{1}(\widetilde{\sigma^{*}L_{\chi}})$ is
    \begin{equation*}
        (r-s)c_{1}(\mathcal{F}) - sc_{1}(\mathcal{E}/\mathcal{F}) = sr \left( \frac{1}{s}c_{1}(\mathcal{F}) - \frac{1}{r}c_{1}(\mathcal{E}) \right).
    \end{equation*}
    Thus semistability of $\mathcal{G}$ implies slope semistability of $\mathcal{E}$.

    Conversely, we claim that to check semistability of $\mathcal{G}$ it suffices to focus on maximal parabolic subgroups.  Indeed, if $P$ is a parabolic subgroup with a dominant character $\chi$, then $\chi$ can be written as a non-negative linear combination of dominant characters of the maximal parabolic subgroups containing $P$.  Since a maximal parabolic subgroup stabilizes a vector subspace, by running the construction in the previous paragraph in reverse we obtain the reverse implication.
\end{proof}

Recall that a principal $SL(r)$-bundle over $U$ corresponds to a rank $r$ locally free sheaf $\mathcal{E}$ on $U$ plus the choice of an isomorphism $\det(\mathcal{E}) \xrightarrow{\cong} \mathcal{O}_{U}$.  In this context we again have a correspondence between semistability of principal bundles and semistability of coherent sheaves.

\begin{lemma} \label{lemm:slvequivalentss}
    Let $X$ be a normal projective variety over an algebraically closed field of characteristic $p$ and let $\alpha \in \Nef^{n-1}(X)$ be a nef curve class.  Let $\mathcal{E}$ be a non-zero torsion-free sheaf of rank $r$ on $X$ and let $U \subset X$ be a big open set contained in the smooth locus of $X$ such that $\mathcal{E}|_{U}$ is locally free.
    
    Assume that $c_{1}(\mathcal{E}) = 0$.  Let $\mathcal{S}$ be any principal $SL(r)$-bundle on $U$ such that $\mathcal{S} \times^{SL(r)} GL(r) \cong \underline{\mathrm{Isom}}_{U}(\mathcal{O}_{U}^{\oplus r},\mathcal{E})$.  Then $\mathcal{E}$ is semistable if and only if $\mathcal{S}$ is semistable.
\end{lemma}

\begin{proof}
    The function that assigns to a parabolic subgroup $P \subset GL(r)$ the parabolic subgroup $P \cap SL(r) \subset SL(r)$ is a bijection that preserves maximality.  For each $P$ we have a surjection of characters $\mathrm{Hom}(P,\mathbb{G}_{m}) \to \mathrm{Hom}(P \cap SL(r),\mathbb{G}_{m})$.  Furthermore, after tensoring by $\mathbb{Q}$ this map induces an isomorphism of the subgroups of characters which vanish along the radical.
    
    Next we translate these properties to principal bundles.  First, there is a bijection between reductions of $\mathcal{S}$ from $SL(r)$ to $P \cap SL(r)$ and reductions of $\mathcal{S} \times^{SL(r)} GL(r)$ from $GL(r)$ to $P$.  Second, the dominant characters of $P \cap SL(r)$ are identified with the dominant characters of $P$ up to rescaling.  With these facts in hand, the proof is essentially the same as before.
\end{proof}

\subsubsection{Preservation of strong semistability}
We start with a quick reminder about inner twists.

\begin{definition} \label{defi:conjaction}
Let $G$ be a connected reductive group and suppose $\mathcal{G} \to U$ is a principal $G$-bundle.  Write $G^{conj}$ for $G$ equipped with the conjugation action.  Then the inner twist is
\begin{equation*}
    \mathcal{G}^{conj} := \mathcal{G} \times^{G} G^{conj}.
\end{equation*}
Note that $\mathcal{G}^{conj}$ is a group scheme over $U$.  Indeed, there is an identification of each fiber of $\mathcal{G}^{conj} \to U$ with $G^{conj}$ using the definition $\mathcal{G}^{conj} = (\mathcal{G} \times G^{conj})/G$.
Furthermore, reductivity is fppf-local on the base, so the fact that $G$ is reductive implies that $\mathcal{G}^{conj}$ is also reductive.

Suppose that $G$ has a left-action on a scheme $T$.  Then $G$ acts on $G \times T$ via the conjugation action on the first factor and the given action on the second factor; with this structure, the action map $G \times T \to T$ is $G$-equivariant.  Thus the $U$-group scheme $\mathcal{G}^{conj}$ acts on the $U$-scheme $\mathcal{G} \times^{G} T$. 
\end{definition}

The key result in this section is the following.

\begin{theorem}[{\cite[Theorem 3.23 and subsequent discussion]{RR84}}] \label{theo:rrtheo}
    Let $X$ be a normal projective variety over an algebraically closed field of characteristic $p$ and let $\alpha \in \Nef^{n-1}(X)$ be a nef curve class.  Let $G$ be a connected reductive group and let $\mathcal{G}$ be a principal $G$-bundle over a big open subset $U \subset X$ contained in the smooth locus.
    
    Suppose $\rho: G \to H$ is a homomorphism of connected reductive groups such that $\rho(Z_{0}(G)) \subset Z_{0}(H)$.  If $\mathcal{G}$ is strongly $\alpha$-semistable, then the principal $H$-bundle $\mathcal{H} := \mathcal{G} \times^{G} H$ is strongly $\alpha$-semistable.
\end{theorem}

For ample polarizations this is proved in \cite[Theorem 3.23]{RR84} and in \cite[Theorem A.3]{Langer04b}.  The proof for arbitrary nef curve classes is exactly the same.

In the proof below we will reference the arguments of \cite[Proposition 3.10.(1)]{RR84}, \cite[Proposition 3.13]{RR84}, and \cite[Last paragraph]{RR84}.  In these propositions, \cite{RR84} establishes that certain line bundles on smooth projective varieties have non-positive degree by exhibiting a section of a tensor power of the dual.  For smooth quasiprojective varieties, the argument goes through unchanged to again exhibit a section of a tensor power of the dual (although of course one cannot define the degree without a projectivity assumption).

\begin{proof}
    Since the construction of $\mathcal{G} \times^{G} H$ commutes with Frobenius pullbacks, it suffices to show that $\mathcal{H}$ is $\alpha$-semistable.

    Suppose that $Q \subset H$ is a parabolic subgroup, that $U' \subset U$ is a big open subset, and that $\sigma: U' \to \mathcal{H}/Q$ is a section corresponding to a reduction of $\mathcal{H}$ from $H$ to $Q$.  To verify Equation \eqref{eq:gss} for all dominant characters $\chi$ on $Q$, by rescaling and taking limits it suffices to consider only those characters such that $\mathcal{T} := L_{\chi}^{-1}$ is a relatively ample line bundle on $\mathcal{H}/Q$ over $U'$. 

    Let $K$ denote the function field of $U'$.  Then $\mathcal{H}/Q_{K}$ is a projective scheme with an ample line bundle $\mathcal{T}_{K}$.  Let $\mathcal{G}^{conj}_{K}$ denote the generic fiber of the $U'$-group scheme $\mathcal{G}^{conj}$.  Since $\mathcal{H}/Q_{K} \cong (\mathcal{G} \times^{G} H/Q)_{K}$, as described in Definition \ref{defi:conjaction} $\mathcal{H}/Q_{K}$ carries an action from the reductive group $\mathcal{G}^{conj}_{K}$. 
    Similarly $\mathcal{G}^{conj}_{K}$ acts on $\mathcal{T}_{K}$ in a compatible way.  The section $\sigma$ defines a rational point $\sigma_{K}$ of $\mathcal{H}/Q_{K}$.
    
    Consider the base change of these constructions to the algebraic closure $\overline{K}$.  The argument splits into two cases, depending on whether or not $\sigma_{\overline{K}}$ is semistable for the $\mathcal{G}^{conj}_{\overline{K}}$-action.
    First suppose that $\sigma_{\overline{K}}$ is semistable.  Then \cite[Proof of Proposition 3.10.(1)]{RR84} shows that for some positive integer $r$ we have $H^{0}(U',\sigma^{*}\mathcal{T}^{\otimes r}) \neq 0$.  Since $U'$ is a big open subset, by passing to the reflexive closure $\widetilde{\sigma^{*}\mathcal{T}}$ we obtain an effective Weil divisor on $X$ representing a multiple of $c_{1}(\widetilde{\sigma^{*}\mathcal{T}})$.  This effective divisor has non-negative intersection against the nef class $\alpha$.  Since $\mathcal{T} = L_{\chi}^{-1}$, this verifies Equation \eqref{eq:gss}.

    Next suppose that $\sigma_{\overline{K}}$ is not semistable.  By \cite{Kempf78} and Galois descent, we can find a destabilizing one parameter subgroup $\lambda$ for $\sigma$ defined over a finite purely inseparable extension $K^{1/p^{e}} \supset K$.  Recall that for the relative Frobenius $F^{e}_{rel}: U'^{(-e)} \to U'$ the function field $K'$ of $U'^{(-e)}$ can be identified as $K' \cong K^{1/p^{e}}$.  Thus for the action of $F^{e*}_{rel}\mathcal{G}^{conj}_{K'}$ on $F^{e*}_{rel}\mathcal{H}/Q_{K'}$, the point $\sigma_{K'}$ has a destabilizing one parameter subgroup $\lambda$ defined over $K'$.  Since $\mathcal{G}$ is strongly semistable, $F^{e*}_{rel}\mathcal{G}$ remains semistable; combined with the assumption that $Z_{0}(G)$ acts trivially on $H/Q$ and $\mathbb{A}^{1}_{\chi}$, we have verified the hypotheses of \cite[Proposition 3.13]{RR84}.  The proof of this result shows that there are positive integers $r,r'$ and auxiliary line bundles $\mathcal{L},\mathcal{M}$ on $U'^{(-e)}$ such that $H^{0}(U'^{(-e)},\mathcal{L} \otimes F^{e*}_{rel}\sigma^{*}\mathcal{T}^{\otimes r}) \neq 0$, $H^{0}(U'^{(-e)},(\mathcal{L}^{\vee})^{\otimes r'} \otimes \mathcal{M}) \neq 0$, and $c_{1}(\widetilde{\mathcal{M}}) \cdot F^{e*}_{rel}\alpha \leq 0$ where $\widetilde{\mathcal{M}}$ is the reflexive extension to $X^{(-e)}$.  By combining these three statements and arguing as before we conclude that the reflexive closure $F^{e\widehat{*}}_{rel}\widetilde{\sigma^{*}\mathcal{T}}$ satisfies $c_{1}(F^{e\widehat{*}}_{rel}\widetilde{\sigma^{*}\mathcal{T}}) \cdot \alpha \geq 0$.   Since the Frobenius rescales intersection numbers by a power of $p$ and $\mathcal{T} = L_{\chi}^{-1}$, this verifies Equation \eqref{eq:gss}.
\end{proof}

We record a number of important corollaries.

\begin{corollary}
    Let $X$ be a normal projective variety over an algebraically closed field of characteristic $p$ and let $\alpha \in \Nef^{n-1}(X)$ be a nef curve class.

    Suppose that $\mathcal{E}$ is a strongly semistable torsion-free sheaf.  Construct the torsion-free sheaf $\mathcal{F}$ by applying the Schur functor $\mathbb{S}_{\lambda}$ associated to the partition $\lambda$ along the locally free locus of $\mathcal{E}$ and taking a reflexive closure.  Then $\mathcal{F}$ is also strongly semistable.
\end{corollary}

In particular, the previous corollary applies to symmetric powers, exterior powers, etc.

\begin{proof}
An element $t \in Z_{0}(GL(V))$ acts on $\mathbb{S}_{\lambda}V$ by $t^{|\lambda|}$ and thus lies in $Z_{0}(GL(\mathbb{S}_{\lambda}V))$.  Thus we may apply Lemma \ref{lemm:equivalentss} and Theorem \ref{theo:rrtheo}.
\end{proof}

\begin{corollary} \label{coro:strongsstensor}
    Let $X$ be a normal projective variety over an algebraically closed field of characteristic $p$ and let $\alpha \in \Nef^{n-1}(X)$ be a nef curve class.
        Suppose that $\mathcal{E}_{1},\mathcal{E}_{2}$ are strongly semistable torsion-free sheaves.  Then $\mathcal{E}_{1} \widehat{\otimes} \mathcal{E}_{2}$ is also strongly semistable.
\end{corollary}

\begin{proof}
     Consider the torsion-free sheaf $\mathcal{E}_{1} \oplus \mathcal{E}_{2}$.  On a big open subset $U \subset X$ contained in the smooth locus, this sheaf defines a principal $GL(V_{1}) \times GL(V_{2})$-bundle $\mathcal{G}$ where $V_{1},V_{2}$ correspond to the two summands.  Recall that every parabolic subgroup of $G := GL(V_{1}) \times GL(V_{2})$ is a product of parabolic subgroups in the two factors.  Thus the strong semistability of $\mathcal{E}_{1}$ and $\mathcal{E}_{2}$ implies the strong semistability of $\mathcal{G}$.  Note that the tensor representation $\rho: G \to GL(V_{1} \otimes V_{2})$ takes the radical $\mathbb{G}_{m} \times \mathbb{G}_{m}$ to the radical $\mathbb{G}_{m}$.  By Theorem \ref{theo:rrtheo}, we conclude that $\mathcal{G} \times^{G} GL(V_{1} \otimes V_{2})$ is strongly semistable.  An application of Lemma \ref{lemm:equivalentss} finishes the proof.
\end{proof}

By combining Theorem \ref{theo:asymsss} and Corollary \ref{coro:strongsstensor} we obtain:

\begin{corollary} \label{coro:lmaxfortensor}
        Let $X$ be a normal projective variety over an algebraically closed field of characteristic $p$ and let $\alpha \in \Nef^{n-1}(X)$ be a nef curve class.

        For any two non-zero torsion-free sheaves $\mathcal{E}_{1},\mathcal{E}_{2}$ we have
        \begin{equation*}
            L^{max}(\mathcal{E}_{1} \widehat{\otimes} \mathcal{E}_{2}) = L^{max}(\mathcal{E}_{1}) + L^{max}(\mathcal{E}_{2}) \qquad \qquad L^{min}(\mathcal{E}_{1} \widehat{\otimes} \mathcal{E}_{2}) = L^{min}(\mathcal{E}_{1}) + L^{min}(\mathcal{E}_{2})
        \end{equation*}
\end{corollary}

\subsubsection{Low height representations} It is also possible to show preservation of semistability at the cost of working in large characteristic.  This perspective was initiated by \cite{Serre94} and developed in our setting by \cite{IMP03} and \cite{BH04}.

\begin{definition}
    Let $G$ be a connected reductive group.  Let $\rho$ be an action of $G$ on a finite-dimensional vector space $V$ and let
    $0 = V_{0} \subset \ldots \subset V_{b} = V$ be a composition series for the $G$-action.  Fix a choice of a maximal torus $T$ and a Borel subgroup $B$ with $T \subset B$ and let $\lambda_{i}$ be the highest weight for the $T$-action on $V_{i}/V_{i-1}$.
    The height of $\rho$ is $\mathrm{ht}(\rho) := 2 \max_{i} \mathrm{ht}(\lambda_{i})$.  Note that $\mathrm{ht}(\rho)$ is independent of the choice of $T \subset B$ and the specific composition series.

    We say that $\rho$ has low height if $\mathrm{ht}(\rho)$ is less than the characteristic $p$ of the ground field.
\end{definition}

Some useful examples of heights of representations are:
\begin{itemize}
    \item For a vector space $V$ of dimension $r$ the $i$th exterior power representation of $SL(V)$ has height $i(r-i)$.
    \item For a vector space $V$ of dimension $r$ the $i$th symmetric power representation of $SL(V)$ has height $i(r-1)$.
    \item Given vector spaces $V_{1},V_{2}$ of dimensions $r_{1},r_{2}$ the tensor product representation of $SL(V_{1}) \times SL(V_{2})$ has height $r_{1}+r_{2} -2$.
\end{itemize}

\begin{theorem}[{\cite[Theorem 3.1]{IMP03}}] \label{theo:lowheightslv}
    Let $X$ be a normal projective variety over an algebraically closed field of characteristic $p$ and let $\alpha \in \Nef^{n-1}(X)$ be a nef curve class.  Let $G$ be a connected semisimple group and let $\mathcal{G}$ be a principal $G$-bundle over a big open subset $U \subset X$ contained in the smooth locus.

    Suppose $\rho: G \to SL(V)$ is a low height representation.  If $\mathcal{G}$ is $\alpha$-semistable, then $\mathcal{G} \times^{G} SL(V)$ is $\alpha$-semistable as a principal $SL(V)$-bundle.
\end{theorem}

\begin{proof}
    Since $Z_{0}(G)$ and $Z_{0}(SL(V))$ are trivial, we are in the setting of Theorem \ref{theo:rrtheo}.  Retaining the notation and argument of the proof, it suffices to verify that when $\sigma_{\overline{K}}$ is not semistable the destabilizing one parameter subgroup $\lambda$ is defined over $K$.  This is proved by \cite[Proof of Theorem 3.1]{IMP03}.
\end{proof}

There are two ways to apply Theorem \ref{theo:lowheightslv} to coherent sheaves.  The first option is to use representation theory to extend Theorem \ref{theo:lowheightslv} to principal $GL(V)$-bundles.  This approach is taken in \cite[Proposition 6.9]{BH04} and the analogous results hold for stability conditions defined by nef curve classes.

The second option is to apply the Bloch-Gieseker construction (\cite{BG71}) to construct principal $SL(r)$-bundles on covers of $X$ from coherent sheaves on $X$.  This method is more suitable for us as it yields slightly better bounds on the characteristic.

\begin{corollary} \label{coro:sslargechar}
        Let $X$ be a normal projective variety over an algebraically closed field of characteristic $p$ and let $\alpha \in \Nef^{n-1}(X)$ be a nef curve class.

        Let $\mathcal{E}_{1},\mathcal{E}_{2}$ be semistable torsion-free sheaves on $X$ of ranks $r_{1},r_{2}$ respectively.
        \begin{enumerate}
            \item If $p > i(r_{1}-i)$ then $\widehat{\Lambda}^{i} \mathcal{E}_{1}$ is semistable.
            \item If $p > i(r_{1}-1)$ then $\widehat{\mathrm{Sym}}^{i} \mathcal{E}_{1}$ is semistable.
            \item If $p > r_{1}+r_{2}-2$ then $\mathcal{E}_{1} \widehat{\otimes} \mathcal{E}_{2}$ is semistable.
        \end{enumerate}
\end{corollary}

\begin{proof}
    (1) Note that if $i \in \{1,r_{1}-1,r_{1}\}$ then the statement is clear, so we may assume otherwise.  In particular, our bound on $p$ implies that $p > r_{1}$.  Let $U$ be the big open locus where $X$ is smooth and $\mathcal{E}_{1}$ is locally free.  Since $\gcd(p,r_{1})=1$, we may apply \cite[Lemma 4.1]{BRS22} to $U$ to obtain a finite separable morphism $f: U' \to U$ such that $f^{*}\det(\mathcal{E}_{1}|_{U})$ is divisible by $r_{1}$ in $\mathrm{Pic}(U')$.  We may take normal projective closures to obtain a finite separable morphism $f: X' \to X$ such that $U'$ is a big open subset of $X'$.  Thus $f^{\widehat{*}}\det(\mathcal{E}_{1})$ is divisible by $r_{1}$ in the class group of $X'$.
    
    By Theorem \ref{theo:separablepullbackisss}, the reflexive pullback of $\mathcal{E}_{1}$ to $X'$ remains semistable.  Choose a Weil divisor $D$ on $X'$ such that $r_{1}D \sim c_{1}(f^{\widehat{*}}\mathcal{E}_{1})$.  Then $f^{\widehat{*}}\mathcal{E}_{1} \widehat{\otimes} \mathcal{O}_{X'}(-D)$ has trivial first Chern class and thus corresponds to a principal $SL(V)$-bundle $\mathcal{S}$ (with $\dim(V) = r_{1}$) on a big open subset of $X'$.  By Lemma \ref{lemm:slvequivalentss} $\mathcal{S}$ is semistable.  Since the exterior power representation has low height, Theorem \ref{theo:lowheightslv} shows that $\mathcal{S} \times^{SL(V)} SL(\Lambda^{i}V)$ is also semistable.  By Lemma \ref{lemm:slvequivalentss} we see that $\widehat{\Lambda}^{i}(f^{\widehat{*}}\mathcal{E}_{1} \widehat{\otimes} \mathcal{O}_{X'}(-D))$ is semistable.  Twisting and applying Theorem \ref{theo:separablepullbackisss} again, we conclude that $\widehat{\Lambda}^{i}\mathcal{E}_{1}$ is semistable.

    (2) and (3) are proved in the analogous way.
\end{proof}

\subsection{Frobenius pullbacks and stability}
We next turn to the study of how slope stability changes under Frobenius pullback.  The general theme is that the behavior of slope stability under Frobenius pullback is controlled by the positivity of the tangent bundle of $X$.

\subsubsection{Cartier correspondence}

The Cartier correspondence describes the behavior of coherent sheaves with respect to the relative Frobenius morphism $F_{rel}: Y \to Y^{(1)}$.  For any $Y^{(1)}$-module $\mathcal{E}$, we define the canonical connection $\nabla_{can, \mathcal{E}}: F_{rel}^{*}\mathcal{E} \to \Omega_{Y} \otimes F_{rel}^{*}\mathcal{E}$ by identifying $F_{rel}^{*}\mathcal{E} = \mathcal{O}_{Y} \otimes_{F_{rel}^{-1}\mathcal{O}_{Y^{(1)}}} \mathcal{E}$ and defining $\nabla_{can,\mathcal{E}}(f \otimes e) = df \otimes e$.

\begin{theorem}[{\cite[Theorem 5.1]{Katz70}}]
    Let $Y$ be a smooth variety over an algebraically closed field of characteristic $p$.  The functor $\mathcal{E} \mapsto (F_{rel}^{*}\mathcal{E},\nabla_{can,\mathcal{E}})$ is an equivalence of categories between the categories of quasicoherent sheaves on $Y^{(1)}$ and the full subcategory of quasicoherent sheaves on $Y$ with an integrable $\mathbb{K}$-connection of zero $p$-curvature.
\end{theorem}

Since $Y^{(1)}$ and $Y$ are abstractly isomorphic, we can identify pullbacks of sheaves under $F_{rel}$ with pullbacks of sheaves under $F_{abs}$.  We will implicitly make this identification in the following statement.

\begin{lemma} \label{lemm:nonzeronabla}
    Let $X$ be a normal projective variety over an algebraically closed field of characteristic $p$.  Let $\mathcal{E}$ be a non-zero torsion-free sheaf on $X$.  Suppose that $\mathcal{F} \subset F^{\widehat{*}}_{abs}\mathcal{E}$ is a saturated torsion-free subsheaf such that there is no torsion-free subsheaf $\mathcal{G} \subset \mathcal{E}$ satisfying $F_{abs}^{\widehat{*}}\mathcal{G} = \mathcal{F}$.
    
    Let $U \subset X$ be a big open subset contained in the smooth locus such that $\mathcal{E}|_{U}$ and $\mathcal{F}|_{U}$ are locally free.  Then $\nabla_{can}$ induces a non-zero $\mathcal{O}_{U}$-linear homomorphism $\mathcal{F}|_{U} \to \Omega_{U} \otimes (F^{*}_{abs}\mathcal{E}/\mathcal{F})|_{U}$.
\end{lemma}

\begin{proof}
    By the Cartier correspondence, we conclude that the restriction $\nabla_{can,\mathcal{E}}|_{\mathcal{F}}: \mathcal{F}|_{U} \to \Omega_{U} \otimes F_{abs}^{*}\mathcal{E}|_{U}$ cannot have image in $\Omega_{U} \otimes \mathcal{F}|_{U}$.  In other words, the induced $\mathcal{O}_{U}$-linear map $\mathcal{F}|_{U} \to \Omega_{U} \otimes (F^{*}_{abs}\mathcal{E}/\mathcal{F})|_{U}$ is non-zero.  
\end{proof}

We will frequently combine the previous result with a stability computation.  (The following statement is most relevant when $e=1$ but will prove useful for arbitrary $e$.)

\begin{proposition}[{\cite[Proposition 1]{SB98}, \cite[Lemma 2.3]{Langer04}}] \label{prop:ssandnonzeronabla}
    Let $X$ be a normal projective variety over an algebraically closed field of characteristic $p$ and let $\alpha \in \Nef^{n-1}(X)$ be a nef curve class.

    Suppose that $\mathcal{E}$ is a semistable torsion-free sheaf such that $F_{abs}^{e\widehat{*}}\mathcal{E}$ is not semistable.  Write the Harder-Narasimhan filtration as
    \begin{equation*}
        0 = \mathcal{F}_{0} \subset \mathcal{F}_{1} \subset \ldots \subset \mathcal{F}_{t} = F^{e\widehat{*}}_{abs}\mathcal{E}.
    \end{equation*}
    Then for each index $i$ with $0 < i < t$ there is no torsion-free sheaf $\mathcal{G} \subset \mathcal{E}$ such that $F^{e\widehat{*}}_{abs}\mathcal{G} = \mathcal{F}_{i}$.
\end{proposition}

\begin{proof}
    Suppose for a contradiction that there were such a subsheaf $\mathcal{G}$.
    Note that $F^{e\widehat{*}}_{abs}$ multiplies the slopes of torsion-free sheaves on $X$ by $p^{e}$.  Then by semistability of $\mathcal{E}$ we must have $\mu_{\alpha}(F^{e\widehat{*}}_{abs}\mathcal{G}) \leq \mu_{\alpha}(F^{e\widehat{*}}_{abs}\mathcal{E})$.  On the other hand, for $0<i<t$ we have $\mu_{\alpha}(\mathcal{F}_{i}) > \mu_{\alpha}(F^{e\widehat{*}}_{abs}\mathcal{E})$ by Lemma \ref{lemm:exactsequence}.  This yields a contradiction.
\end{proof}

\subsubsection{Slope inequalities}

The following result was first proved in a weaker form by \cite[Corollary 6.2]{Langer04}.  The argument below gives a tighter bound following ideas of \cite[Proposition 3.9]{Sun10} and \cite[Theorem 2.7]{Trivedi12}.

\begin{theorem} \label{theo:mainfrobpullback}
    Let $X$ be a normal projective variety over an algebraically closed field of characteristic $p$ and let $\alpha \in \Nef^{n-1}(X)$ be a nef curve class.

    Suppose that $\mathcal{E}$ is a non-zero torsion-free sheaf on $X$.  Choose a positive integer $e$ so that each successive quotient in the Harder-Narasimhan filtration of $F_{abs}^{e\widehat{*}}\mathcal{E}$ is strongly semistable.  Write this filtration as
    \begin{equation*}
        0 = \mathcal{F}_{0,\infty} \subset \mathcal{F}_{1,\infty} \subset \ldots \subset \mathcal{F}_{\ell,\infty} = F^{e\widehat{*}}_{abs}\mathcal{E}.
    \end{equation*}
    Let $K \subset \{1,\ldots,\ell-1\}$ be the indices of the terms in this sequence such that there is a torsion-free subsheaf $\mathcal{G} \subset \mathcal{E}$ satisfying $F^{e\widehat{*}}_{abs}(\mathcal{G}) = \mathcal{F}_{i,\infty}$ and write $s = |K|$.
    \begin{enumerate}
        \item For any $k \in \{1,\ldots,\ell-1 \} \backslash K$, we have
        \begin{equation*}
            \frac{1}{p^{e}} \left(\mu_{\alpha}(\mathcal{F}_{k,\infty}/\mathcal{F}_{k-1,\infty}) - \mu_{\alpha}(\mathcal{F}_{k+1,\infty}/\mathcal{F}_{k,\infty}) \right) \leq \frac{1}{p} \max\{ \mu_{\alpha}^{max}(\widehat{\Omega_{X}}), 0\}.
        \end{equation*}
        \item We have
    \begin{equation*}
        L^{max}(\mathcal{E}) - L^{min}(\mathcal{E}) \leq \frac{(\ell -1)- s}{p} \max \{ \mu_{\alpha}^{max}(\widehat{\Omega_{X}}), 0 \} + (\mu^{max}_{\alpha}(\mathcal{E}) - \mu_{\alpha}^{min}(\mathcal{E})).
    \end{equation*}
    \end{enumerate}
\end{theorem}

\begin{proof}
    (1) Fix an index $k \not \in K$. For a given $e' < e$, there may or may not be a saturated subsheaf $\mathcal{G} \subset F_{abs}^{e'\widehat{*}}\mathcal{E}$ such that $F_{abs}^{(e-e')\widehat{*}}\mathcal{G} = \mathcal{F}_{k,\infty}$.  Let $e_{k} < e$ be the largest non-negative integer such that there is \emph{not} such a sheaf $\mathcal{G}$.  This means that there \emph{is} such a sheaf $\mathcal{G}$ for $e_{k}+1$.
    
    Let $U \subset X$ be a big open set contained in the smooth locus such that $\mathcal{G}|_{U}$ and $\mathcal{E}|_{U}$ are locally free.  Then Lemma \ref{lemm:nonzeronabla} shows that we have a non-zero $\mathcal{O}_{U}$-morphism
    \begin{equation*}
        \mathcal{G}|_{U} \to \Omega_{U} \otimes F^{(e_{k}+1)\widehat{*}}_{abs}\mathcal{E}/\mathcal{G}|_{U}.
    \end{equation*}
    Applying adjunction and taking reflexive closures, we see that we also have a non-zero morphism
    \begin{equation*}
        \mathcal{G} \widehat{\otimes} (F^{(e_{k}+1)\widehat{*}}_{abs}\mathcal{E}/\mathcal{G})^{\vee} \to \widehat{\Omega_{X}}.
    \end{equation*}
    By combining Lemma \ref{lemm:slopeinequalityzerohom} with the comparison $L^{min} \leq \mu^{min}$, we obtain an inequality
    \begin{equation*}
        L^{min}(\mathcal{G} \widehat{\otimes} (F^{(e_{k}+1)\widehat{*}}_{abs}\mathcal{E}/\mathcal{G})^{\vee}) \leq \mu_{\alpha}^{max}(\widehat{\Omega_{X}}).
    \end{equation*}
    By combining Lemma \ref{lemm:dual} and Corollary \ref{coro:lmaxfortensor} we can reinterpret the left-hand side to obtain
    \begin{equation*}
        L^{min}(\mathcal{G}) - L^{max}(F^{(e_{k}+1)\widehat{*}}_{abs}\mathcal{E}/\mathcal{G}) \leq \mu_{\alpha}^{max}(\widehat{\Omega_{X}}).
    \end{equation*}     
    By construction, we know that
    \begin{align*}
        L^{min}(\mathcal{G}) &= \frac{1}{p^{e-e_{k}-1}} L^{min}(\mathcal{F}_{k,\infty}) = \frac{1}{p^{e-e_{k}-1}}\mu_{\alpha}(\mathcal{F}_{k,\infty}/\mathcal{F}_{k-1,\infty}) \\
        L^{max}(F^{(e_{k}+1)\widehat{*}}_{abs}\mathcal{E}/\mathcal{G}) & = \frac{1}{p^{e-e_{k}-1}} L^{max}(F^{e\widehat{*}}_{abs}\mathcal{E}/\mathcal{F}_{k,\infty}) = \frac{1}{p^{e-e_{k}-1}} \mu_{\alpha}(\mathcal{F}_{k+1,\infty}/\mathcal{F}_{k,\infty})
    \end{align*}
    Finally we conclude
    \begin{equation*}
        \frac{1}{p^{e}}\mu_{\alpha}(\mathcal{F}_{k,\infty}/\mathcal{F}_{k-1,\infty}) - \frac{1}{p^{e}} \mu_{\alpha}(\mathcal{F}_{k+1,\infty}/\mathcal{F}_{k,\infty}) \leq \frac{1}{p^{e_{k}+1}} \mu^{max}_{\alpha}(\widehat{\Omega_{X}}).
    \end{equation*}
    Since $e_{k} \geq 0$ for every $k \not \in K$, we obtain the desired statement.

    (2) We next analyze $\mathcal{F}_{k,\infty}$ for $k \in K$ in the case when $K$ is non-empty.  Write $K = \{k_{1},\ldots,k_{s} \}$ where the indices are in increasing order and set $k_{0} = 0$ and $k_{s+1} = \ell$.  For each index $k_{i} \in K$, there is a saturated subsheaf $\mathcal{G}_{i} \subset \mathcal{E}$ so that $\mathcal{F}_{k_{i},\infty} = F^{e\widehat{*}}_{abs}\mathcal{G}_{i}$.  Set $\mathcal{G}_{0} = 0$ and $\mathcal{G}_{s+1} = \mathcal{E}$.  Since the Frobenius morphism is faithfully flat on the smooth locus, we have inclusions $\mathcal{G}_{i} \subset \mathcal{G}_{i+1}$ for every $i$.  Then
    \begin{align*}
        \mu_{\alpha}(\mathcal{F}_{k_{i},\infty}/\mathcal{F}_{k_{i}-1,\infty}) & = \mu^{min}_{\alpha}(\mathcal{F}_{k_{i},\infty}/\mathcal{F}_{k_{(i-1)},\infty})  \\
        & \leq \mu_{\alpha}(\mathcal{F}_{k_{i},\infty}/\mathcal{F}_{k_{(i-1)},\infty}) \\
        & = p^{e} \mu_{\alpha}(\mathcal{G}_{i}/\mathcal{G}_{i-1}).
    \end{align*}
    Similarly $\mu_{\alpha}(\mathcal{F}_{k_{i}+1,\infty}/\mathcal{F}_{k_{i},\infty}) \geq p^{e} \mu_{\alpha}(\mathcal{G}_{i+1}/\mathcal{G}_{i})$.  Thus
    \begin{align*}
        \sum_{i \in K} \frac{1}{p^{e}} (\mu_{\alpha}(\mathcal{F}_{i,\infty}/\mathcal{F}_{i-1,\infty}) - \mu_{\alpha}(\mathcal{F}_{i+1,\infty}/\mathcal{F}_{i,\infty})) & \leq \sum_{i=1}^{s} \mu_{\alpha}(\mathcal{G}_{i}/\mathcal{G}_{i-1}) - \mu_{\alpha}(\mathcal{G}_{i+1}/\mathcal{G}_{i}) \\
        & \leq \mu_{\alpha}^{max}(\mathcal{E}) - \mu_{\alpha}^{min}(\mathcal{E}).
    \end{align*}
    To estimate $L^{max}(\mathcal{E}) - L^{min}(\mathcal{E})$, we sum up the differences in slopes between successive quotients occurring in the filtration $\{\mathcal{F}_{k,\infty}\}_{k=0}^{\ell}$.  We separate out this sum into the terms in $K$ and the terms not in $K$ and use (1) to estimate the latter contribution.  Thus we have:
    \begin{align*}
        L^{max}(\mathcal{E}) - L^{min}(\mathcal{E}) & = \sum_{i \not \in K} \frac{1}{p^{e}} (\mu_{\alpha}(\mathcal{F}_{i,\infty}/\mathcal{F}_{i-1,\infty}) - \mu_{\alpha}(\mathcal{F}_{i+1,\infty}/\mathcal{F}_{i,\infty})) \\
        & \qquad + \sum_{i \in K} \frac{1}{p^{e}} (\mu(\mathcal{F}_{i,\infty}/\mathcal{F}_{i-1,\infty}) - \mu(\mathcal{F}_{i+1,\infty}/\mathcal{F}_{i,\infty})) \\
        & \leq \frac{(\ell-1)-s}{p} \max \{ \mu_{\alpha}^{max}(\widehat{\Omega_{X}}), 0 \} + (\mu_{\alpha}^{max}(\mathcal{E}) - \mu_{\alpha}^{min}(\mathcal{E})).
    \end{align*}
\end{proof}

Since for every positive integer $e$ we have $L^{max}(\mathcal{E}) \geq \mu_{\alpha}^{max}(F_{abs}^{e\widehat{*}}\mathcal{E})/p^{e}$ and $L^{min}(\mathcal{E}) \leq \mu_{\alpha}^{min}(F_{abs}^{e\widehat{*}}\mathcal{E})/p^{e}$, the previous theorem also gives bounds on the instability of $\mathcal{E}$ under any iterated Frobenius pullback.  Thus the previous result subsumes many different statements in the literature.  Below are a couple examples.

\begin{corollary}[{\cite[Corollary 6.2]{Langer04}}] \label{coro:lmaxmin}
    Let $X$ be a normal projective variety over an algebraically closed field of characteristic $p$ and let $\alpha \in \Nef^{n-1}(X)$ be a nef curve class.

    Suppose that $\mathcal{E}$ is a semistable sheaf on $X$.  Then
    \begin{equation*}
        L^{max}(\mathcal{E}) - L^{min}(\mathcal{E}) \leq \frac{\rk(\mathcal{E})-1}{p} \max \{ \mu_{\alpha}^{max}(\widehat{\Omega_{X}}), 0 \}.
    \end{equation*}
\end{corollary}

\begin{proof}
    Let $e$ be a sufficiently large integer so that successive quotients in the Harder-Narasimhan filtration of $F^{e\widehat{*}}_{abs}\mathcal{E}$ are strongly semistable.  Write this filtration as
    \begin{equation*}
        0 = \mathcal{F}_{0,\infty} \subset \mathcal{F}_{1,\infty} \subset \ldots \subset \mathcal{F}_{\ell,\infty} = F^{e\widehat{*}}_{abs}\mathcal{E}.
    \end{equation*}
    By Proposition \ref{prop:ssandnonzeronabla}, for $0 < k < \ell$ there is no torsion-free subsheaf $\mathcal{G} \subset \mathcal{E}$ such that $\mathcal{F}_{k,\infty} = F_{abs}^{e\widehat{*}}\mathcal{G}$.  Thus Theorem \ref{theo:mainfrobpullback}.(2) shows that
        \begin{equation*}
        L^{max}(\mathcal{E}) - L^{min}(\mathcal{E}) \leq \max \left\{ \frac{\ell-1}{p}\mu_{\alpha}^{max}(\widehat{\Omega_{X}}),0 \right\}.
    \end{equation*}
    We conclude by noting that $\ell \leq \rk(\mathcal{E})$.
\end{proof}

\begin{corollary}[{\cite[Theorem 2.1]{MR83}, \cite[Corollary 3.10]{Sun10}}]
    Let $X$ be a normal projective variety over an algebraically closed field of characteristic $p$ and let $\alpha \in \Nef^{n-1}(X)$ be a nef curve class.
    
    Suppose that $\mu_{\alpha}^{min}(T_{X}) \geq 0$.  Then every semistable sheaf $\mathcal{E}$ is strongly semistable.
\end{corollary}

\begin{proof}
    By Corollary \ref{coro:lmaxmin} we see that $L^{max}(\mathcal{E}) = L^{min}(\mathcal{E})$.  This implies strong semistability.
\end{proof}

We can also avoid working with iterated Frobenius pullbacks at the cost of assuming the characteristic is large.

\begin{proposition} \label{prop:strongerbound}
Let $X$ be a normal projective variety over an algebraically closed field of characteristic $p$ and let $\alpha \in \Nef^{n-1}(X)$ be a nef curve class.

Let $\mathcal{E}$ be a semistable torsion-free sheaf on $X$ such that $p > \rk(\mathcal{E})-2$.  Suppose that $F_{abs}^{\widehat{*}}\mathcal{E}$ is not semistable and write
\begin{equation*}
    0 = \mathcal{F}_{0} \subset \mathcal{F}_{1} \subset \ldots \subset \mathcal{F}_{t} = F_{abs}^{\widehat{*}}\mathcal{E}
\end{equation*}
for the Harder-Narasimhan filtration of $F_{abs}^{\widehat{*}}\mathcal{E}$. Then for every $1 \leq i < t$ we have
\begin{equation*}
    \mu_{\alpha}(\mathcal{F}_{i}/\mathcal{F}_{i-1}) - \mu_{\alpha}(\mathcal{F}_{i+1}/\mathcal{F}_{i}) \leq \mu_{\alpha}^{max}(\widehat{\Omega_{X}}).
\end{equation*}
\end{proposition}

\begin{proof}
    Fix an index $i$ with $0 < i < t$.  Let $U \subset X$ be a big open subset contained in the smooth locus such that $\mathcal{E}|_{U}$ and $\mathcal{F}_{i}|_{U}$ are locally free.  By combining Lemma \ref{lemm:nonzeronabla} and Proposition \ref{prop:ssandnonzeronabla} we have a non-zero morphism
    \begin{equation*}
        \mathcal{F}_{i}|_{U} \to \Omega_{U} \otimes F^{*}_{abs}\mathcal{E}/\mathcal{F}_{i}|_{U}.
    \end{equation*}
    Applying adjunction and taking reflexive closures, we see that we also have a non-zero morphism
    \begin{equation*}
        \mathcal{F}_{i} \widehat{\otimes} (F^{\widehat{*}}_{abs}\mathcal{E}/\mathcal{F}_{i})^{\vee} \to \widehat{\Omega_{X}}.
    \end{equation*}    
    Thus we must have
    \begin{equation*}
        \mu_{\alpha}^{min}(\mathcal{F}_{i} \widehat{\otimes} (F^{\widehat{*}}_{abs}\mathcal{E}/\mathcal{F}_{i})^{\vee}) \leq \mu_{\alpha}^{max}(\widehat{\Omega_{X}}).
    \end{equation*}
    By our characteristic assumption, Lemma \ref{lemm:dual} and Corollary \ref{coro:sslargechar} show that the leftmost term is
    \begin{equation*}
        \mu_{\alpha}^{min}(\mathcal{F}_{i}) - \mu_{\alpha}^{max}(F^{\widehat{*}}_{abs}\mathcal{E}/\mathcal{F}_{i}) = \mu_{\alpha}(\mathcal{F}_{i}/\mathcal{F}_{i-1}) - \mu_{\alpha}(\mathcal{F}_{i+1}/\mathcal{F}_{i}).
    \end{equation*}
\end{proof}

\bibliographystyle{alpha}
\bibliography{normalstab}
\end{document}